\documentclass[11pt]{article}
\usepackage[margin=1in]{geometry}
\usepackage[T1]{fontenc}
\usepackage[utf8]{inputenc}
\usepackage{amsmath,amssymb,amsfonts,amsthm,mathtools}
\usepackage{algorithm}
\usepackage{algpseudocode}
\usepackage{microtype}
\usepackage{xcolor}
\usepackage{lmodern}
\usepackage{tikz}
\usepackage{pgfplots}
\usepackage{placeins}
\usepackage{booktabs}
\usepackage{array}
\usepackage{tabularx}
\usepackage{pifont}
\usepackage[numbers,square,sort]{natbib}
\makeatletter
\def\NAT@spacechar{~}
\makeatother
\usepackage[colorlinks=true,linkcolor=blue,citecolor=blue,urlcolor=blue]{hyperref}
\newcommand{\authormail}[1]{{\small\hypersetup{urlcolor=black}\href{mailto:#1}{\texttt{#1}}}}
\usepackage{orcidlink}
\pgfplotsset{compat=1.18}
\hypersetup{pdftitle={Limited-Memory First-Order Methods: Identification, Reconstruction, $H$-Duality and Application to Optimized Methods},pdfauthor={Pierre Vernimmen, Fran\c{c}ois Glineur, Baptiste Goujaud, Adrien B. Taylor}}

\theoremstyle{plain}
\newtheorem{theorem}{Theorem}[section]
\newtheorem{proposition}[theorem]{Proposition}
\newtheorem{lemma}[theorem]{Lemma}
\newtheorem{corollary}[theorem]{Corollary}
\newtheorem{sublemma}{Lemma}[subsection]
\theoremstyle{remark}
\newtheorem{remark}[theorem]{Remark}
\theoremstyle{definition}
\newtheorem{definition}[theorem]{Definition}

\algrenewcommand\algorithmicindent{1.35em}

\newcommand{\R}{\mathbb R}
\newcommand{\ip}[2]{\left\langle #1,#2\right\rangle}
\newcommand{\norm}[1]{\left\lVert #1\right\rVert}
\newcommand{\Fcl}{\mathcal F_{\mu,L}}
\newcommand{\rsum}{r_{\mathrm{sum}}}
\newcommand{\Hd}{H^{\mathrm A}}
\DeclareMathOperator{\rank}{rank}
\DeclareMathOperator{\Diag}{Diag}
\DeclareMathOperator{\dist}{dist}
\newcommand{\one}{\mathbf 1}
\newcommand{\sym}{\operatorname{sym}}
\newcommand{\cmark}{\textcolor{green!55!black}{\ding{51}}}
\newcommand{\xmark}{\textcolor{red!75!black}{\ding{55}}}
\newcommand{\backmatter}{}

\begin{document}
	\raggedbottom
\nocite{bansalgupta2019,barretaylorbach2023,bokaltschuler2026,dasguptavanparysryu2024,dewildeveen1998,diakonikolaswang2022,drori2017exact,drori2026weighted,droritaylor2020efficient,droritaylor2022oracle,droriteboulle2014,eidelmangohberg1999,eidelmangohberghaimovici2014,fazlyabribeiromoraripreciado2018,goujaud2024pepit,grimmerjopark2026,kaashoekwoerdeman1989,kamrihendrickxglineur2025,kim2021accelerated,kimfessler2016ogm,kimfessler2021ogmg,kimryudasgupta2026,kim2023timereversed,kornowski2024open,leeparkryu2021,lessardrechtpackard2016,lessardseiler2020,lessardudell2026,nemirovskyyudin1983,nesterov1983,nesterov2004,odonoghuecandes2015,rodmanwoerdeman1995,schererebenbauer2021,shuwang2026unified,taylordrori2023item,taylorhendrickxglineur2017,taylorvanscoylessard2018,upadhyaya2026proxitem,vanderveendewilde1994lowcomplexity,vanscoy2018tmm,vandebrilvanbarelmastronardi2005,vandebrilvanbarelmastronardi2008,vernimmen2026anytime,weibel2026optimized,woerdeman1987,yoongrimmer2026,yoonkimsuhryu2024,yoonryugrimmer2025,zhang2025anytime}

\title{Limited-Memory First-Order Methods: Identification, Reconstruction, $H$-Duality and Application to Optimized Methods}

\author{
\makebox[0.42\textwidth]{Pierre Vernimmen}\\
{\normalsize ICTEAM Institute, UCLouvain}\\
\authormail{pierre.vernimmen@uclouvain.be}
\and
\makebox[0.42\textwidth]{Fran\c{c}ois Glineur}\\
{\normalsize ICTEAM Institute, UCLouvain}\\
{\normalsize LIDAM/CORE, UCLouvain}\\
\authormail{francois.glineur@uclouvain.be}
\and
\makebox[0.42\textwidth]{Baptiste Goujaud}\\
{\normalsize T\'el\'ecom SudParis}\\
{\normalsize Institut Polytechnique de Paris}\\
\authormail{baptiste.goujaud@telecom-sudparis.eu}
\and
\makebox[0.42\textwidth]{Adrien B. Taylor}\\
{\normalsize Universit\'e Grenoble Alpes, Inria, CNRS}\\
{\normalsize Grenoble INP, LJK}\\
\authormail{adrien.taylor@inria.fr}}

\date{}

\maketitle

\begin{abstract}
In general, linear-span first-order methods add at every iteration a combination of all observed gradients to the previous iterate. However, many methods of interest can be implemented in a way that requires storing only a few auxiliary vectors, instead of all past gradients. We call such methods limited-memory first-order methods, and these auxiliary vectors correspond to sums of past gradients. 

A given fixed-coefficient first-order method can be described by a lower triangular coefficient matrix, which expresses each step as a combination of previous gradients. We show that the minimum number of running sums required to implement a method is equal to the minimum rank of a completion of the strict lower triangular part of its coefficient matrix, and propose a reconstruction algorithm. Moreover, taking the $H$-dual of a method, which reflects the coefficient matrix across its anti-diagonal, preserves its memory and translates into a permutation and time-reversal on the coefficients of the running sum recurrences.

We apply our results to two optimized first-order methods for smooth strongly convex functions. First, our reconstruction applied to a method numerically designed by Shu and Wang reveals a recurrence with one auxiliary vector, which is related to the Triple Momentum Method (TMM) of Van Scoy, Freeman and Lynch with a particular initialization. We prove that the reconstructed method, called TMM-GD, is minimax optimal at every horizon for the final distance to the minimizer relative to the initial gradient norm, within all linear-span first-order methods including adaptive ones.  

Second, the Information-Theoretic Exact Method (ITEM) of Taylor and Drori, optimal for the final distance relative to the initial one, is known to be memory one. We give an explicit recurrence with one auxiliary vector for its $H$-dual, called ITEM-GG, and prove that ITEM-GG is the unique minimax optimal nonadaptive linear-span method for the final gradient norm relative to the initial one.


\end{abstract}

\newpage
{\hypersetup{linkcolor=black}\tableofcontents}
\newpage

\section{Introduction}

When analyzing or designing first-order methods, it is useful to describe each iterate using a linear combination of all the gradients computed so far. This makes the coefficients explicit and allows comparing methods with different structures in a common framework. However, implementing this representation directly requires storing the entire gradient history, which quickly becomes impractical as the number of iterations grows. The same method may still admit a short recurrence using a small number of auxiliary vectors. Nesterov's method \citep{nesterov1983}, the Optimized Gradient Method (OGM) \citep{kimfessler2016ogm}, the Triple Momentum Method (TMM) \citep{vanscoy2018tmm} and the Information-Theoretic Exact Method (ITEM) \citep{taylordrori2023item}, for instance, run with only one auxiliary vector, although their coefficient matrices, introduced below, have a (nearly) dense lower triangle. We study how to identify and recover such limited-memory implementations from the general structure.

Fix a horizon of $N$ iterations. We consider nonadaptive linear-span first-order methods, where a method is said to be \emph{nonadaptive} if its coefficients do not depend on any oracle information. For a differentiable function $f:\mathbb R^d\to\mathbb R$, such a method can be written as
\begin{equation}\label{eq:fsfom}
x_{k+1}=x_k-\frac{1}{\bar L}\sum_{i=0}^{k}H_{k,i}\nabla f(x_i),
\qquad k=0,\ldots,N-1,
\end{equation}
where $\bar L>0$ is any reference scale.\footnote{We write $\bar L$ rather than $L$ because Sections~\ref{sec:memory} and~\ref{sec:recon} involve no smoothness constant. Only $H/\bar L$ enters~\eqref{eq:fsfom}, so choosing $\bar L$ amounts to choosing the scale of $H$. From Section~\ref{sec:tmmgd} on, $\bar L=L$.} The matrix $H\in\mathbb R^{N\times N}$ is lower triangular because the step producing $x_{k+1}$ can only use gradients evaluated at $x_0,\dots,x_k$. Its entries may depend on the horizon $N$. Our structural results concern $H$, so they remain valid when the gradients are replaced by arbitrary vectors, e.g., the outputs of a proximal oracle.

The implementations we look for store a few auxiliary vectors, which they update linearly with each new gradient. We call a method that can be run with $r$ such vectors a \emph{memory-$r$ method}. When the auxiliary vectors are running sums of past gradients, we show, using a result from \citet{woerdeman1987}, that the smallest number of these sums is equal to the minimum rank of a completion of the strict lower triangle of $H$.

Two situations in particular call for starting from $H$. First, $H$ is the natural object of analysis and design. It is standard in Performance Estimation (PEP), which expresses the worst-case performance of an optimization method as an optimization problem \citep{droriteboulle2014,taylorhendrickxglineur2017}, and methods are designed by minimizing this worst-case over $H$, a nonconvex problem in general. Such design is as old as PEP: \citet{droriteboulle2014} solve a convex relaxation of this problem numerically, and \citet{kimfessler2016ogm} solve the relaxation analytically, which gives OGM. Convex relaxations have also given, for example, ITEM \citep{taylordrori2023item}, an inexact proximal point method \citep{barretaylorbach2023} and projection-free methods for online learning \citep{weibel2026optimized}, while numerical solutions of the nonconvex problem have given OGM-G \citep{kimfessler2021ogmg} and an accelerated proximal point method for maximally monotone operators \citep{kim2021accelerated}. Dedicated algorithms now solve the nonconvex problem globally \citep{dasguptavanparysryu2024} or locally \citep{kamrihendrickxglineur2025}. The designed method is then a dense matrix, from which a practical implementation must be recovered, and nothing in the design favors short recurrences. Yet many of these methods turn out a posteriori to admit one, and these implementations are derived case by case, as in recent works \citep{bokaltschuler2026,kimryudasgupta2026}. \citet{droritaylor2020efficient} also observe that the coefficients of the methods they construct numerically have a rank-one structure.

Second, $H$-duality is defined on $H$: it maps a method to the method whose coefficient matrix is the anti-transpose of $H$, that is, $H$ reflected across its anti-diagonal. \citet{kimfessler2021ogmg} show that the coefficient matrices of OGM and OGM-G are related in this way, and \citet{kim2023timereversed} introduce $H$-duality for all methods of the form~\eqref{eq:fsfom}. It relates methods that reduce the function value of a smooth convex function to methods that reduce its gradient norm, and \citet{shuwang2026unified} develop it into a unified theory for smooth strongly convex minimization and contractive fixed-point problems.

The two uses meet in a method that \citet{shuwang2026unified} design numerically for $L$-smooth, $\mu$-strongly convex functions under the gradient-to-distance criterion, i.e., the worst case of $\|x_N-x_\star\|^2/\|\nabla f(x_0)\|^2$, where $x_N$ is the output and $x_\star$ is the minimizer. For five iterations, they report its coefficient matrix rounded to four decimals (Section~\ref{sec:decoding}). From these rounded values, our reconstruction procedure recovers a method that uses only one auxiliary vector.

The contributions of the paper are as follows.
\begin{enumerate}
\item \textbf{Memory and $H$-duality of general recurrences (Section~\ref{sec:memory}).} We characterize when general recurrences need fewer stored vectors than running sums, and we observe that the two counts coincide for almost every coefficient matrix and for all methods of this paper, so that the simpler running sums are generally enough. We then write the $H$-dual directly from the recurrence, without ever expanding the matrix $H$, and we show that it is also a memory-$r$ method.

\item \textbf{Limited-memory reconstruction and implementation (Section~\ref{sec:recon}).} We adapt the minimum-rank completion of \citet{woerdeman1987} into an algorithm that reconstructs a method from its coefficient matrix. For exact coefficients, it returns an implementation with the fewest running sums. For rounded coefficients, a tolerance lets it ignore residuals at the scale of the rounding, and we give an exact formula for the resulting coefficient error. We implement methods in a momentum form, which avoids the potentially badly scaled coefficients of running sums, and we show that the $H$-dual of a momentum form is again a momentum form, with its coefficients in reverse order.

\item \textbf{Discovery and minimax optimality of TMM-GD (Section~\ref{sec:tmmgd}).} Shu and Wang see their numerical design as a path to an analytic description of a method that they call ITEM-GD. We derive the closed form of the method that our reconstruction reveals: TMM with a particular initialization, which we call TMM-GD. We prove that, on the class $\Fcl$ of $L$-smooth, $\mu$-strongly convex functions, it is minimax optimal for the gradient-to-distance criterion over all deterministic, possibly adaptive, linear-span methods, in dimension at least $N+1$. Its recurrence is time-invariant, and its output rule does not depend on $N$, so a single run attains the minimax value at every horizon.

\item \textbf{$H$-duality applied to ITEM: ITEM-GG and TMM-GG (Section~\ref{sec:item}).} \citet{shuwang2026unified} define ITEM-GG implicitly as the $H$-dual of ITEM through its coefficient matrix and derive its gradient-norm guarantee from a PEP certificate of ITEM. With the $H$-dual momentum form of Section~\ref{sec:duality-momentum}, we write ITEM-GG explicitly as a memory-one recurrence and show its guarantee with a potential (Lyapunov) argument. We prove that, on $\Fcl$, ITEM-GG is the unique minimax optimal method for the gradient-to-gradient criterion $\norm{\nabla f(x_N)}^2/\norm{\nabla f(x_0)}^2$ among nonadaptive linear-span methods~\eqref{eq:fsfom}, in dimension at least $N+1$, which proves for these methods the conjectured optimality of its rate \citep{shuwang2026unified}. Since its first step differs for any two horizons, minimax optimality at a horizon $N\ge2$ requires knowing $N$ before the first query. As $N\to\infty$, its coefficients tend to those of a stationary method, TMM-GG, which is also the $H$-dual of TMM with zero initial momentum. We prove for TMM-GG an accelerated guarantee at every horizon of a single run, and a better bound for a second output formed from the same state.
\end{enumerate}

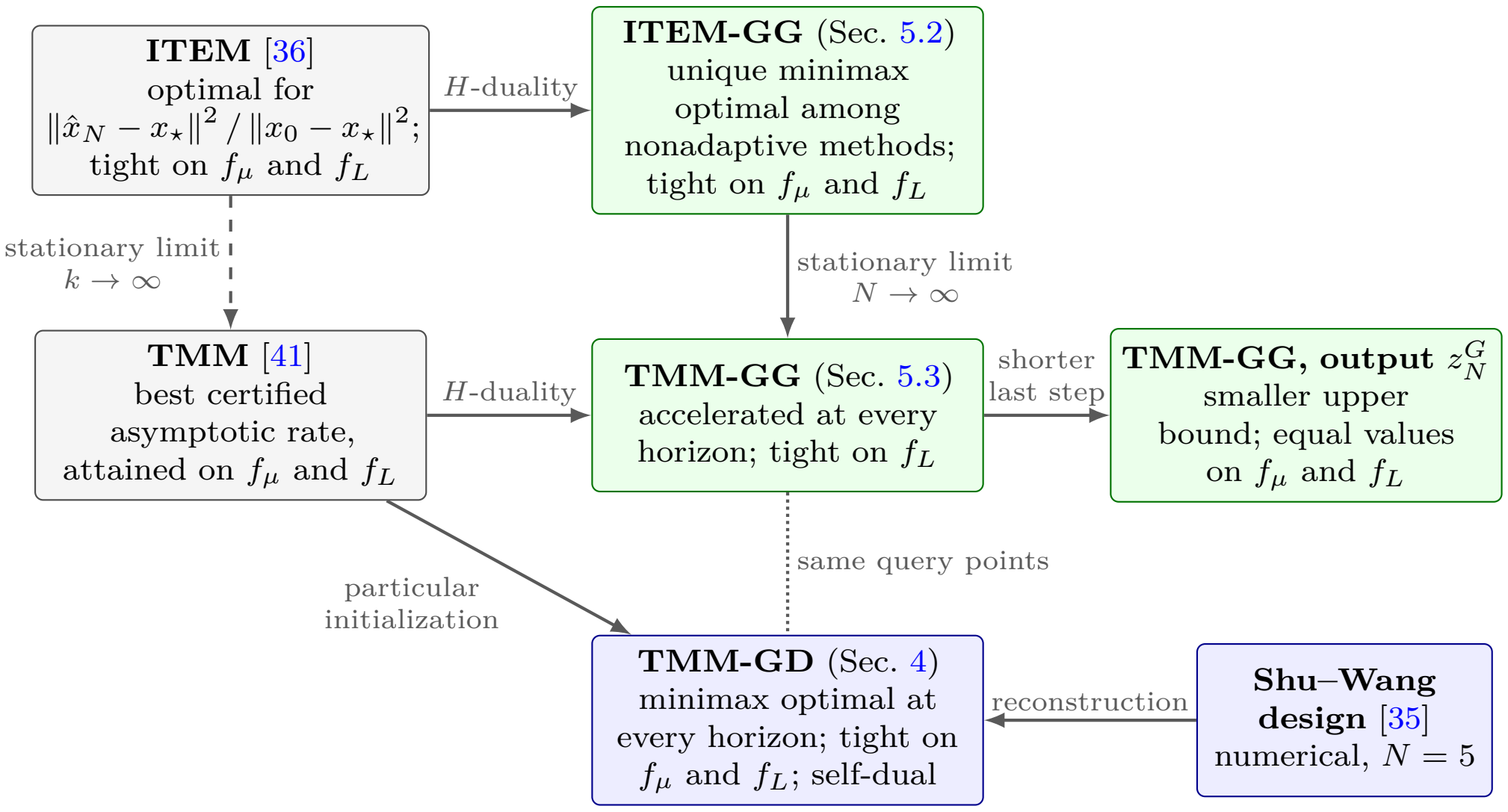
\begin{figure}[H]
\centering
\begin{tikzpicture}[
  x=0.85cm,font=\fontsize{8}{9.5}\selectfont,
  box/.style={draw,rounded corners=2pt,align=center,text width=3.25cm,inner sep=3pt,minimum height=1.35cm,execute at begin node={\hyphenpenalty=10000\exhyphenpenalty=10000}},
  gd/.style={box,fill=blue!7,draw=blue!55!black},
  gg/.style={box,fill=green!8,draw=green!45!black},
  other/.style={box,fill=gray!8,draw=gray!60!black},
  arr/.style={-latex,thick,gray!70!black},
  lab/.style={font=\fontsize{7}{8}\selectfont,align=center,inner sep=2pt},
]
\node[other,text width=3.3cm] (item) at (0,0) {\textbf{ITEM} \cite{taylordrori2023item}\\ optimal for $\norm{\hat x_N-x_\star}^2/\norm{x_0-x_\star}^2$; tight on $f_\mu$ and $f_L$};
\node[gg]    (itemgg) at (5.8,0)      {\textbf{ITEM-GG} (Sec.~\ref{sec:itemgg})\\ unique minimax optimal among nonadaptive methods; tight on $f_\mu$ and $f_L$};
\node[other] (tmm)    at (0,-2.7)     {\textbf{TMM} \cite{vanscoy2018tmm}\\ best certified asymptotic rate, attained on $f_\mu$ and $f_L$};
\node[gg]    (tmmgg)  at (5.8,-2.7)   {\textbf{TMM-GG} (Sec.~\ref{sec:tmmgg})\\ accelerated at every horizon; tight on $f_L$};
\node[gg]    (zg)     at (11.2,-2.7)  {\textbf{TMM-GG, output $z^G_N$}\\ smaller upper bound; equal values on $f_\mu$ and $f_L$};
\node[gd]    (tmmgd)  at (5.8,-5.4)   {\textbf{TMM-GD} (Sec.~\ref{sec:tmmgd})\\ minimax optimal at every horizon; tight on $f_\mu$ and $f_L$; self-dual};
\node[gd,text width=2.4cm] (sw) at (11.6,-5.4) {\textbf{Shu--Wang design} \cite{shuwang2026unified}\\ numerical, $N=5$};
\draw[arr] (item) -- node[lab,above] {$H$-duality} (itemgg);
\draw[arr,dashed] (item) -- node[lab,left] {stationary limit\\ $k\to\infty$} (tmm);
\draw[arr] (itemgg) -- node[lab,right] {stationary limit\\ $N\to\infty$} (tmmgg);
\draw[arr] (tmm) -- node[lab,above] {$H$-duality} (tmmgg);
\draw[arr] (tmmgg) -- node[lab,above] {shorter\\ last step} (zg);
\draw[arr] (tmm) -- node[lab,below left] {particular\\ initialization} (tmmgd);
\draw[thick,gray!70!black,densely dotted] (tmmgg) -- node[lab,right] {same query points} (tmmgd);
\draw[arr] (sw) -- node[lab,above] {reconstruction} (tmmgd);
\end{tikzpicture}
\caption{The methods of Sections~\ref{sec:tmmgd} and~\ref{sec:item} and their links. Grey (suffix DD): distance-to-distance criterion $\norm{\hat x_N-x_\star}^2/\norm{x_0-x_\star}^2$; blue (suffix GD): gradient-to-distance criterion $\norm{\hat x_N-x_\star}^2/\norm{\nabla f(x_0)}^2$; green (suffix GG): gradient-to-gradient criterion $\norm{\nabla f(\hat x_N)}^2/\norm{\nabla f(x_0)}^2$, where $\hat x_N$ is the output. The quadratics $f_\mu:=\frac\mu2\norm{\cdot-x_\star}^2$ and $f_L:=\frac L2\norm{\cdot-x_\star}^2$ are the extreme functions of the class $\Fcl$.}
\label{fig:map}
\end{figure}

Although Section~\ref{sec:tmmgd} starts from a numerical coefficient matrix and Section~\ref{sec:item} from ITEM, the methods TMM-GD and TMM-GG share the same query points and only differ in their output, and the potential of TMM-GD certifies both (Section~\ref{sec:tmmgg}). Figure~\ref{fig:map} summarizes these methods and their relations, and Table~\ref{tab:summary} summarizes their guarantees.

\begin{table}[H]
\caption{Worst-case guarantees on $\Fcl$, valid for every $N\ge1$, with $q=\mu/L$, $\rho=1-\sqrt q$ and $A_N$ given by~\eqref{eq:item-A}. The column $\hat x_N$ gives the output at which the criterion is evaluated: the last iterate $x_N$, or one of the points $z_N$ and $z^G_N$ of~\eqref{eq:output} and~\eqref{eq:output-gg}. For TMM-GD, optimal means minimax optimal among possibly adaptive linear-span methods. For ITEM-GG, it means minimax optimal among nonadaptive linear-span methods, and no other such method is optimal. These two statements hold in dimension $d\ge N+1$. The last column says whether the query points are independent of $N$.}
\label{tab:summary}
\centering
\setlength{\tabcolsep}{2.5pt}%
\def\rs{\ensuremath{\vcenter{\hbox{\rule{0pt}{27pt}}}}}
\def\down#1{\raisebox{-13.5pt}[0pt][0pt]{#1}}
\def\two#1#2{\begin{tabular}[c]{@{}c@{}}#1\\#2\end{tabular}}
\def\thickrule{\noalign{\hrule height 0.8pt}}%
\def\dotrule{\omit&\multispan{7}\leaders\hbox to 1.5pt{\hss\vrule height 0.4pt width 0.4pt\hss}\hfill\cr\noalign{\vskip-0.4pt}}
\begin{tabularx}{\textwidth}{@{}c|c|c|c|>{\raggedright\arraybackslash}c|c|c|c@{}}
\thickrule
\ensuremath{\vcenter{\hbox{\rule{0pt}{22pt}}}}Criterion & Method & $\hat x_N$ & Bound & Status & Source & \begin{tabular}[c]{@{}c@{}}Stationary\\interior\end{tabular} & \begin{tabular}[c]{@{}c@{}}Independent\\of $N$\end{tabular}\\
\thickrule
\down{$\dfrac{\norm{\hat x_N-x_\star}^2}{\norm{x_0-x_\star}^2}$} & ITEM & $x_N$ & \rs$\dfrac1{1+qA_N}$ & optimal & \citep{droritaylor2022oracle,taylordrori2023item} & \xmark & \cmark\\
\dotrule
& TMM & $z_N$ & \rs$\dfrac{\rho^{2N}}{q}$ & tight & \cite{vanscoy2018tmm} & \cmark & \cmark\\
\thickrule
$\dfrac{\norm{\hat x_N-x_\star}^2}{\norm{\nabla f(x_0)}^2}$ & TMM-GD & $z_N$ & \rs$\dfrac{(1+\sqrt q)^2\rho^{2N}}{4\mu^2}$ & optimal & Thm.~\ref{thm:minimax} & \cmark & \cmark\\
\thickrule
 & ITEM-GG & $x_N$ & \rs$\dfrac1{1+qA_N}$ & optimal & \cite{shuwang2026unified}, Thm.~\ref{thm:itemgg-minimax} & \xmark & \xmark\\
\dotrule
$\dfrac{\norm{\nabla f(\hat x_N)}^2}{\norm{\nabla f(x_0)}^2}$ & TMM-GG & $x_N$ & \rs$\dfrac{\rho^{2N}}{q}$ & tight & Thm.~\ref{thm:tmmgg} & \cmark & \cmark\\
\dotrule
 & TMM-GG & $z^G_N$ & \rs$\dfrac{(1+\sqrt q)^2\rho^{2N}}{4q}$ & \two{upper}{bound} & Thm.~\ref{thm:tmmgg-out} & \cmark & \cmark\\
\thickrule
\end{tabularx}
\end{table}

This work draws on minimum-rank completion \citep{woerdeman1987,rodmanwoerdeman1995} and realization theory, which constructs minimal representations of structured matrices and time-varying systems \citep{kaashoekwoerdeman1989,dewildeveen1998,eidelmangohberghaimovici2014} and approximates them with low-complexity representations \citep{vanderveendewilde1994lowcomplexity}. For time-invariant optimization algorithms, \citet{lessardudell2026} use state-space realizations and transfer functions to characterize equivalence between algorithms and to construct minimal realizations. Sections~\ref{sec:memory} and~\ref{sec:recon} build on minimum-rank completion and realization theory for methods with a finite horizon, whose coefficients may vary with the iteration.

This work also relates to control theory, where a method is a linear system in feedback with the gradient, and integral quadratic constraints (IQC) or Lyapunov functions certify its performance \citep{lessardrechtpackard2016,fazlyabribeiromoraripreciado2018,taylorvanscoylessard2018,bansalgupta2019}. Among time-invariant methods, TMM \citep{vanscoy2018tmm} attains the best asymptotic rate certified in this way \citep{lessardseiler2020,schererebenbauer2021}.

The appendices contain the remaining proofs.

\section{Memory and \texorpdfstring{$H$}{H}-duality of general recurrences}
\label{sec:memory}

In this section and the next, $H\in\R^{N\times N}$ is any lower-triangular matrix, and the gradients in~\eqref{eq:fsfom} are replaced by arbitrary vectors $u_0,\dots,u_{N-1}$ of $\R^d$:
\begin{equation}\label{eq:Hmap}
x_{k+1}=x_k-\frac{1}{\bar L}\sum_{i=0}^{k}H_{k,i}u_i,
\qquad k=0,\dots,N-1.
\end{equation}

\subsection{General recurrences and running sums}
\label{sec:memory-completion}
 
\begin{definition}[general recurrence and running-sum form]\label{def:running-sum}
A \emph{general recurrence} of \eqref{eq:Hmap} with $r$ stored vectors writes it as
\begin{equation}\label{eq:general}
x_{k+1}=x_k-\frac{1}{\bar L}\bigl(S_kc_k+d_ku_k\bigr),
\qquad
S_{k+1}=S_kA_k+u_kb_k^\top,
\qquad S_0=0,
\end{equation}
for some $A_k\in\R^{r\times r}$, $b_k,c_k\in\R^r$ and $d_k\in\R$, which may depend on $N$. The $r$ columns of $S_k\in\R^{d\times r}$ are its \emph{stored vectors}. A \emph{running-sum form} is a general recurrence with $A_k=I$ for every $k$:
\begin{equation}\label{eq:additive}
x_{k+1}=x_k-\frac{1}{\bar L}\bigl(S_kc_k+d_ku_k\bigr),
\qquad
S_{k+1}=S_k+u_kb_k^\top,
\qquad S_0=0.
\end{equation}
\end{definition}
 
At step $k$, the method reads the stored vectors with the weights $c_k$, and then replaces them by the combinations $S_kA_k$ and adds $u_k$ to them with the weights $b_k$. A general recurrence is thus a linear time-varying system with state $S_k$, input $u_k$ and output $S_kc_k+d_ku_k$, whose input--output map is $H$, and the memory defined below is the smallest size of its state. The matrices $A_k$ are unrelated to the scalars $A_k$ of ITEM in Section~\ref{sec:item}.
 
Expanding $S_k$ in \eqref{eq:general} gives
\begin{equation}\label{eq:general-factor}
H_{k,i}=b_i^\top A_{i+1}\cdots A_{k-1}c_k\quad(i<k),
\qquad
H_{k,k}=d_k,
\end{equation}
where the product is the identity when $i=k-1$. Conversely, any factorization~\eqref{eq:general-factor} of $H$ with $A_k\in\R^{r\times r}$ defines a general recurrence of~\eqref{eq:Hmap} with $r$ stored vectors. In the terminology of \citet{eidelmangohberg1999}, the vectors $c_k^\top$, the matrices $A_k^\top$ and the vectors $b_k$ are lower quasiseparable generators of order $r$ for $H$. We use this correspondence to write the $H$-dual of a recurrence (Proposition~\ref{prop:dual-additive}) and to compute the memory from the blocks of $H$ (Theorem~\ref{thm:addmemory}).
 
\begin{definition}[memory]\label{def:memory}
The method \eqref{eq:Hmap} is a \emph{memory-$r$ method} if it admits a general recurrence with $r$ stored vectors. The \emph{memory} $r(H)$ of $H$ is the smallest such $r$, and its \emph{running-sum memory} $\rsum(H)$ is the smallest number of stored vectors of a running-sum form.
\end{definition}
 
We do not count the iterate $x_k$, the current oracle output $u_k$, or the coefficients in $r$. When every transition $A_k$ is invertible, a general recurrence can be rewritten as a running-sum form with the same number of stored vectors: with $P_0=I$ and $P_{k+1}=P_kA_k$, the matrices $\hat S_k:=S_kP_k^{-1}$ satisfy $\hat S_0=0$ and, since $P_{k+1}^{-1}=A_k^{-1}P_k^{-1}$, $\hat S_{k+1}=\hat S_k+u_k\hat b_k^\top$ with $\hat b_k:=P_{k+1}^{-\top}b_k$, while $S_kc_k=\hat S_k\hat c_k$ with $\hat c_k:=P_kc_k$.

\subsection{\texorpdfstring{$H$}{H}-duality}
\label{sec:duality}

\begin{definition}[$H$-dual, \cite{kim2023timereversed}]\label{def:hdual}
The \emph{anti-transpose} of $H$ is the matrix $\Hd$ with entries $\Hd_{k,i}=H_{N-1-i,\,N-1-k}$. The \emph{$H$-dual} of \eqref{eq:Hmap} is the method with coefficient matrix $\Hd$.
\end{definition}

The anti-transpose reflects $H$ across its anti-diagonal, maps lower-triangular matrices to lower-triangular matrices, and satisfies $(\Hd)^{\mathrm A}=H$.

\begin{proposition}[$H$-dual of a general recurrence]\label{prop:dual-additive}
If \eqref{eq:general} is a general recurrence of \eqref{eq:Hmap} with $r$ stored vectors, then the $H$-dual of \eqref{eq:Hmap} has the general recurrence
\begin{equation}\label{eq:dual-additive}
x_{k+1}=x_k-\frac{1}{\bar L}\bigl(\tilde S_kb_{N-1-k}+d_{N-1-k}u_k\bigr),
\quad
\tilde S_{k+1}=\tilde S_kA_{N-1-k}^\top+u_kc_{N-1-k}^\top,
\quad \tilde S_0=0,
\end{equation}
with $r$ stored vectors, which is a running-sum form if \eqref{eq:general} is one. In particular, $r(\Hd)=r(H)$ and $\rsum(\Hd)=\rsum(H)$.
\end{proposition}

\begin{proof}
The form \eqref{eq:dual-additive} has the coefficients $\tilde b_k=c_{N-1-k}$, $\tilde c_k=b_{N-1-k}$ and $\tilde d_k=d_{N-1-k}$, and the transitions $\tilde A_k=A_{N-1-k}^\top$. Let $i<k$, so that $N-1-i>N-1-k$. By \eqref{eq:general-factor}, transposed,
\[
\Hd_{k,i}=H_{N-1-i,\,N-1-k}=c_{N-1-i}^\top A_{N-2-i}^\top\cdots A_{N-k}^\top b_{N-1-k}=\tilde b_i^\top\tilde A_{i+1}\cdots\tilde A_{k-1}\tilde c_k,
\]
and $\Hd_{k,k}=d_{N-1-k}=\tilde d_k$. Hence \eqref{eq:dual-additive} runs the $H$-dual, and $r(\Hd)\le r(H)$. Applying this to $\Hd$, whose anti-transpose is $H$, gives $r(H)\le r(\Hd)$. If $A_k=I$ for every $k$, then $\tilde A_k=I$, and the same argument gives $\rsum(\Hd)=\rsum(H)$.
\end{proof}

At step $k$, the $H$-dual stores $u_k$ with the weights $c_{N-1-k}$ with which the original method reads its stored vectors at step $N-1-k$, and it reads with the weights $b_{N-1-k}$ with which the original method stores $u_{N-1-k}$. $H$-duality thus reverses time and exchanges storing and reading: \eqref{eq:dual-additive} is the transposed system of \eqref{eq:general}, run backward in time.

\subsection{Computing the memories}
\label{sec:memory-blocks}

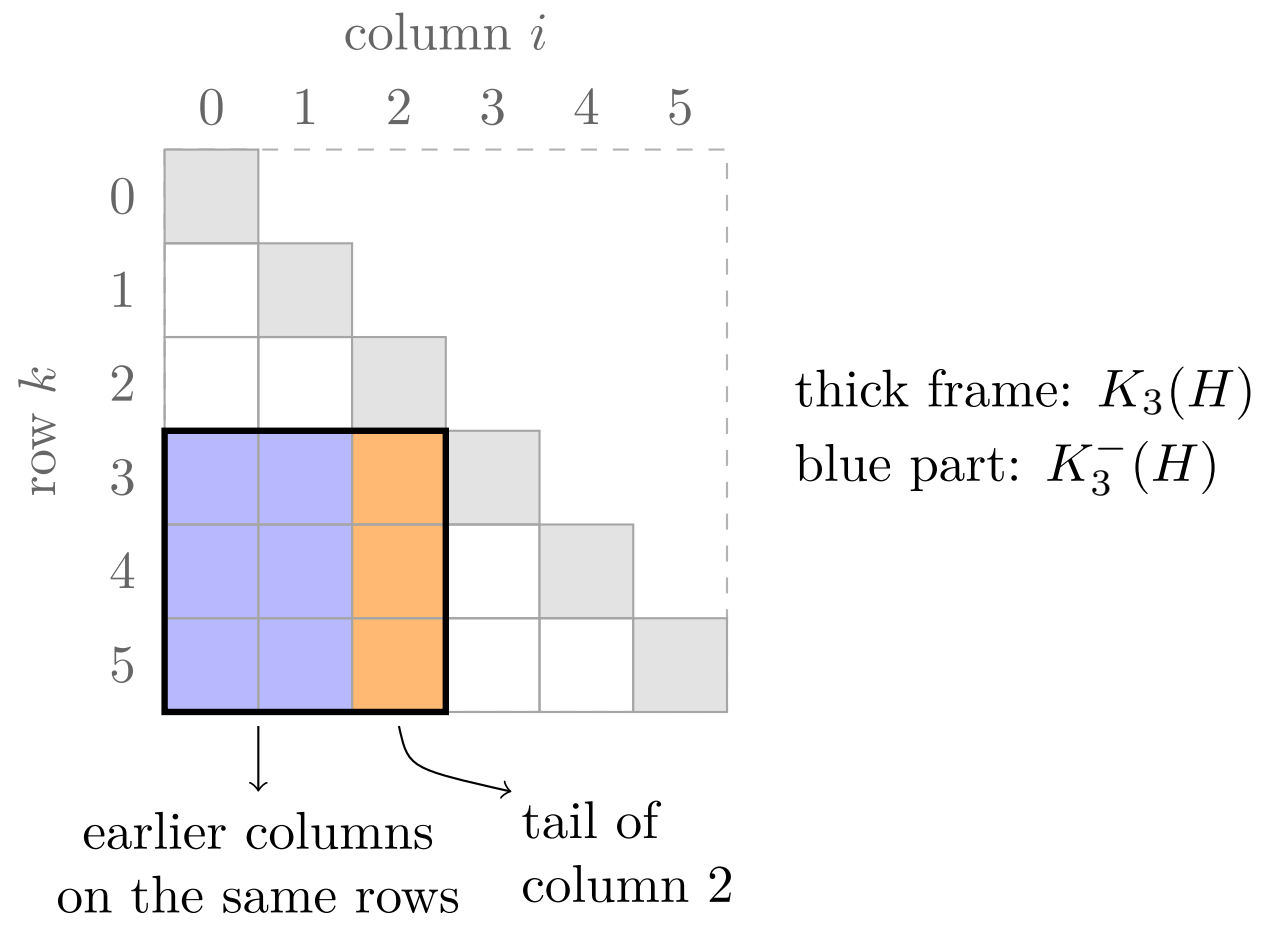
\begin{figure}[h]
\centering
\begin{tikzpicture}[x=0.62cm,y=0.62cm,font=\small]
  \draw[gray!60,dashed] (0,0) rectangle (6,-6);
  \foreach \k in {0,...,5}{\foreach \i in {0,...,\k}{
    \ifnum\i=\k \fill[gray!22] (\i,-\k) rectangle ++(1,-1); \fi
    \draw[gray!70] (\i,-\k) rectangle ++(1,-1);}}
  \fill[blue!28] (0,-3) rectangle (2,-6);
  \fill[orange!55] (2,-3) rectangle (3,-6);
  \foreach \k in {3,...,5}{\foreach \i in {0,...,2}{\draw[gray!70] (\i,-\k) rectangle ++(1,-1);}}
  \draw[very thick] (0,-3) rectangle (3,-6);
  \foreach \i in {0,...,5}{\node[gray!80!black] at (\i+0.5,0.45) {$\i$};}
  \foreach \k in {0,...,5}{\node[gray!80!black] at (-0.45,-\k-0.5) {$\k$};}
  \node[gray!80!black,anchor=south] at (3,0.85) {column $i$};
  \node[gray!80!black,rotate=90,anchor=south] at (-0.95,-3) {row $k$};
  \draw[->,thin] (1,-6.15) -- (1,-6.85) node[below,align=center] {earlier columns\\ on the same rows};
  \draw[->,thin] (2.5,-6.15) .. controls (2.6,-6.6) .. (3.7,-6.85) node[below right=-2pt,align=left] {tail of\\ column $2$};
  \node[anchor=west,align=left] at (6.5,-3.0) {thick frame: $K_3(H)$\\[2pt] blue part: $K_3^-(H)$};
\end{tikzpicture}
\caption{The test~\eqref{eq:mergetest} for column $i=2$ and $N=6$. Column $2$ needs no new stored vector exactly when its tail (orange) is a combination of the earlier columns on the same rows (blue). The diagonal (gray) and the entries above it play no role.}
\label{fig:mergetest}
\end{figure}

To compute $\rsum(H)$, follow the vectors $u_0,u_1,\dots$ in the order in which they arrive. The later steps use $u_i$ with the weights $H_{i+1:N-1,\,i}$, the entries of column $i$ below the diagonal, which we call the \emph{tail} of column $i$. Suppose that this tail is a linear combination of the earlier columns on the same rows, i.e.,
\begin{equation}\label{eq:mergetest}
H_{i+1:N-1,\,i}\in\operatorname{span}\bigl\{H_{i+1:N-1,\,0},\dots,H_{i+1:N-1,\,i-1}\bigr\}.
\end{equation}
Then no new stored vector is needed. Otherwise, column $i$ is \emph{new}: $u_i$ requires a new stored vector. To count these vectors, let
\begin{equation*}
K_k(H):=H_{k:N-1,\,0:k-1},
\qquad
K_k^-(H):=H_{k:N-1,\,0:k-2}
\end{equation*}
be the weights (visualized in Figure~\ref{fig:mergetest}) that steps $k,\dots,N-1$ give to $u_0,\dots,u_{k-1}$, and the same block without its last column (empty for $k=1$, with rank zero). The tail of column $k-1$ is the last column of $K_k(H)$, and the earlier columns on the same rows form $K_k^-(H)$. Column $k-1$ is therefore new exactly when $\rank K_k(H)>\rank K_k^-(H)$, and the running-sum memory is the number of new columns. This is the lower-triangular minimum-rank completion theorem, for which we include a short proof, whose construction motivates our reconstruction algorithm.

The block $K_k(H)$ holds the weights that steps $k,\dots,N-1$ give to $u_0,\dots,u_{k-1}$, so after step $k-1$, a general recurrence must store $\rank K_k(H)$ independent combinations of these vectors. A running sum may need more, since it keeps every combination that it has stored. We say that step $k\ge1$ \emph{forgets} if deleting the first row of $K_k(H)$ lowers its rank. This row holds the weights of step $k$ on $u_0,\dots,u_{k-1}$, so after step $k$, one combination of these vectors is no longer needed. In the next theorem, the maximum and the sum are zero when $N=1$.

\begin{theorem}[memory and running-sum memory, after \cite{kaashoekwoerdeman1989,dewildeveen1998} for (i) and \cite{woerdeman1987} for (ii)]\label{thm:addmemory}
For every lower-triangular $H\in\R^{N\times N}$, \\
(i) $r(H)=\max_{1\le k\le N-1}\rank K_k(H)$; \\
(ii) the running-sum memory is
\begin{equation}\label{eq:radd}
\rsum(H)=\sum_{k=1}^{N-1}\bigl(\rank K_k(H)-\rank K_k^-(H)\bigr);
\end{equation}
(iii) $r(H)\le\rsum(H)$, with equality if and only if $i<k$ whenever column $i$ is new and step $k$ forgets.
\end{theorem}

\begin{proof}
For (ii), build a running-sum form column by column. If the tail of column $i$ is a combination $\sum_{j<i}\omega_jH_{i+1:N-1,\,j}$ of the earlier columns on the same rows, as in~\eqref{eq:mergetest}, the weights $b_i=\sum_{j<i}\omega_jb_j$ reproduce this tail. Otherwise, a new stored vector, read with the weights of this tail, does. This uses one stored vector per positive term of \eqref{eq:radd}.

Conversely, take a running-sum form with $r$ stored vectors. By \eqref{eq:general-factor} with $A_k=I$, $K_k(H)=C_k^\top[\,b_0\ \cdots\ b_{k-1}\,]$ and $K_k^-(H)=C_k^\top[\,b_0\ \cdots\ b_{k-2}\,]$, with $C_k=[\,c_k\ \cdots\ c_{N-1}\,]$. A positive term thus requires $b_{k-1}$ to raise the dimension of $\operatorname{span}\{b_0,\dots,b_{k-1}\}\subseteq\R^r$, which happens for at most $r$ indices.

Parts (i) and (iii) are proved in Appendix~\ref{app:general}.
\end{proof}

Low memory is not generic: for almost every lower-triangular $H$, $r(H)=\rsum(H)=\lfloor N/2\rfloor$ (Corollary~\ref{cor:generic}). Thus, for $N\ge4$, the coefficient matrices of memory at most one form a set of zero Lebesgue measure. Yet many optimal methods have memory one, such as OGM and ITEM, and so do the methods studied in this work.

Running sums never overwrite a stored vector, whereas general recurrences can, through singular transitions $A_k$, and may then need fewer vectors. Since a general recurrence whose transitions are all invertible can be rewritten as a running-sum form with the same number of stored vectors (Section~\ref{sec:memory-completion}), fewer than $\rsum(H)$ stored vectors require a singular transition. By Theorem~\ref{thm:addmemory}(iii), this only helps methods that forget part of the past before new information arrives. Restart schemes forget by design: for a memory-one method restarted on a fixed schedule \citep{odonoghuecandes2015}, running sums need one stored vector per segment, whereas a general recurrence reuses a single one, with $A_k=0$ and $c_k=0$ at the restarts. For methods that do not forget in this way, i.e., for almost every coefficient matrix, the two counts coincide, so we work from now on with running sums, which are simpler to characterize and to reconstruct.

\section{From coefficient matrices to limited-memory implementations}
\label{sec:recon}

\subsection{Reconstruction}
\label{sec:reconstruction}

For a running-sum form, $S_k=\sum_{i<k}u_ib_i^\top$, so comparing \eqref{eq:additive} with \eqref{eq:Hmap} gives
\begin{equation}\label{eq:addfactor}
H_{k,i}=c_k^\top b_i\quad(i<k),
\qquad
H_{k,k}=d_k.
\end{equation}
The diagonal of $H$ is matched by $d_k$, whatever $b_k$ and $c_k$ are, so only the strict lower triangle constrains these coefficients. Collect them in $B=[\,b_0\ \cdots\ b_{N-1}\,]$ and $C=[\,c_0\ \cdots\ c_{N-1}\,]$ in $\R^{r\times N}$. By \eqref{eq:addfactor}, $C^\top B$ agrees with $H$ below the diagonal and has rank at most $r$. Conversely, any matrix of rank at most $r$ that agrees with $H$ below the diagonal admits this factorization and defines a memory-$r$ implementation. Thus, determining $\rsum(H)$ is a minimum-rank completion problem in which the entries on and above the diagonal are free. By Theorem~\ref{thm:addmemory}(ii), its minimum rank is the number of new columns, and the construction in the proof of this theorem gives the factors $B$ and $C$ of a completion with this rank. In the terminology of structured linear algebra, \eqref{eq:addfactor} is a diagonal plus semiseparable representation of $H$, with generators $B$ and $C$ \citep{vandebrilvanbarelmastronardi2005,vandebrilvanbarelmastronardi2008,eidelmangohberghaimovici2014}.

Algorithm~\ref{alg:additive} implements the construction in the proof of Theorem~\ref{thm:addmemory}(ii), with the test~\eqref{eq:mergetest} relaxed by a relative tolerance $\varepsilon\ge0$. Storing $u_i$ with the weights $b$ gives it the future coefficients $C_{:,\,i+1:N-1}^\top b$, which the algorithm fits to the tail of column $i$ by least squares. Column $i$ is \emph{merged} when the relative residual is at most $\varepsilon$ (line~5). Otherwise, the algorithm introduces a new stored vector.

\begin{algorithm}[h!]
\caption{Running-sum reconstruction}
\label{alg:additive}
\begin{algorithmic}[1]
\Require lower-triangular $H$, relative tolerance $\varepsilon\ge0$
\Ensure coefficients $b_k$, $c_k$ and $d_k$ for $k=0,\dots,N-1$, and the matrix $G$ they define
\State $C\gets$ matrix with no rows
\For{$i=0,\dots,N-2$}
\State $\bar b\gets$ the \emph{minimum-norm} minimizer of $\|C_{:,\,i+1:N-1}^\top b-H_{i+1:N-1,\,i}\|_2$ over $b\in\R^r$, where $r$ is the number of rows of $C$
\State $e_i\gets H_{i+1:N-1,\,i}-C_{:,\,i+1:N-1}^\top\bar b$
\If{$\|e_i\|_2\le\varepsilon\,\|H_{i+1:N-1,\,i}\|_2$}
\State $b_i\gets\bar b$
\Else
\State append one zero coordinate to every previously defined $b_j$, $j<i$
\State append the row $(0,\dots,0,H_{i+1,i},\dots,H_{N-1,i})$ to $C$
\State $b_i\gets(0,\dots,0,1)^\top\in\R^{r+1}$
\EndIf
\EndFor
\State $b_{N-1}\gets0$, $c_k\gets C_{:,k}$ and $d_k\gets H_{k,k}$ for $k=0,\dots,N-1$
\State $G_{k,i}\gets c_k^\top b_i$ for $0\le i<k\le N-1$, $G_{k,k}\gets d_k$ for $k=0,\dots,N-1$, and $G_{k,i}\gets0$ for $i>k$
\State \Return $(b_k,c_k,d_k)_{k=0}^{N-1}$ and $G$
\end{algorithmic}
\end{algorithm}

When $C$ has no rows, $\bar b$ is empty and $e_i$ is the whole tail. The returned matrix $G$ does not depend on the minimum-norm choice in line~3. Each least-squares problem has at most $r$ unknowns, where $r$ is the final number of rows of $C$, so the algorithm costs $O((r+1)^2N^2)$ operations, independently of the dimension $d$.

\begin{theorem}[reconstruction guarantees]\label{thm:additive-error}
Let $G$ be the matrix returned by Algorithm~\ref{alg:additive}, and let $e_i$ be the residual computed at each merged column $i$. Then, \\
(i) for $\varepsilon=0$, $G=H$ and $C$ has $\rsum(H)$ rows: the reconstruction is exact and uses the minimum number of stored vectors of a running-sum form; \\
(ii) for every $\varepsilon\ge0$, the running-sum memory of $G$ equals the number of rows of $C$; \\
(iii) for every $\varepsilon\ge0$,
\begin{equation}\label{eq:additive-error}
\|H-G\|_F^2=\sum_{\text{merged columns }i}\|e_i\|_2^2,
\end{equation}
and consequently
\begin{equation}\label{eq:additive-bounds}
\|H-G\|_F\le\varepsilon\,\|H\|_F,
\qquad
\max_{i<k}|H_{k,i}-G_{k,i}|\le\varepsilon\max_{i}\|H_{i+1:N-1,\,i}\|_2.
\end{equation}
\end{theorem}

\begin{proof}
We first prove (iii). If column $i$ is new, its tail is copied into a new row of $C$ and $b_i$ selects this row, so $G_{i+1:N-1,\,i}=H_{i+1:N-1,\,i}$. If column $i$ is merged, then $G_{i+1:N-1,\,i}=C_{:,\,i+1:N-1}^\top b_i$, so $H_{i+1:N-1,\,i}-G_{i+1:N-1,\,i}=e_i$. Rows added to $C$ later do not change this, because $b_i$ receives a zero for each of them. The diagonal is exact by the choice of $d_k$. Summing the squared errors column by column gives \eqref{eq:additive-error}. The test in line~5 gives $\|e_i\|_2\le\varepsilon\|H_{i+1:N-1,\,i}\|_2$ at every merged column, which implies \eqref{eq:additive-bounds}.

For (ii), on rows $i+1,\dots,N-1$, the earlier columns of $G$ span the same space as the rows of $C$ created before step $i$. Hence column $i$ of $G$ satisfies \eqref{eq:mergetest} for $G$ exactly when it was merged, since $e_i\neq0$ for a new column. Theorem~\ref{thm:addmemory}(ii), applied to $G$, then shows that its running-sum memory is the number of rows of $C$.

For (i), the test in line~5 with $\varepsilon=0$ forces $e_i=0$ at every merged column. Hence $G=H$ by (iii), and $C$ has $\rsum(H)$ rows by (ii).
\end{proof}

\begin{remark}[numerical matrices]\label{rem:numerical}
Rounding can destroy exact column dependencies, so a rounded matrix can have a larger memory than the method it approximates. The tolerance $\varepsilon$ lets the reconstruction ignore residuals at the scale of the rounding, at the cost of the coefficient error bounded in Theorem~\ref{thm:additive-error}(iii), so it should exceed the relative size of the rounding errors on the tails of the columns, as in Section~\ref{sec:decoding}. In floating-point arithmetic, the computation adds its own rounding errors, which Theorem~\ref{thm:additive-error} does not account for, so $\varepsilon=0$ is too strict even for exact input, for which we use $\varepsilon=10^{-8}$. For $\varepsilon>0$, the algorithm is greedy, as it decides each column once, in order, so $G$ need not have the smallest running-sum memory among the matrices $G'$ with $\|H-G'\|_F\le\varepsilon\|H\|_F$. This can happen when merging a column requires large combination coefficients, which amplify the perturbation. In the following experiment, this occurred rarely: for running-sum forms with $N\in\{10,20,40\}$, $r\in\{1,2,3,5\}$ stored vectors and standard Gaussian coefficients, whose coefficient matrices we perturbed by independent relative Gaussian errors of standard deviation $\varepsilon/100$ in each entry, with $\varepsilon\in\{10^{-8},10^{-4},10^{-2}\}$, the algorithm returned at most $r$ stored vectors in at least $97\%$ of the $500$ trials of each case, and in all trials for $r=1$. Since anti-transposition preserves the Frobenius norm and the running-sum memory (Proposition~\ref{prop:dual-additive}), one can also apply the algorithm to $\Hd$ and anti-transpose the result, which gives another reconstruction $G$ with $\|H-G\|_F\le\varepsilon\|H\|_F$, and keep the one with fewer stored vectors.
\end{remark}

\subsection{Implementation}
\label{sec:implementation}

Running sums expose the storage structure of a method, but their coefficients can be badly scaled: for a method with momentum, they contain cumulative products of the momentum coefficients (Proposition~\ref{prop:momentum} below), which, for ITEM, grow or decay geometrically, so that a running-sum implementation eventually overflows in floating-point arithmetic. To implement methods, we therefore use the following form, which avoids these products.

\begin{definition}[momentum form]\label{def:momentum}
A \emph{momentum form} of \eqref{eq:Hmap} with $r$ momenta writes it as
\begin{equation}\label{eq:local}
M_{k+1}=\bigl(u_k\one^\top+M_k\Theta_k\bigr)W_k,
\qquad
x_{k+1}=x_k-\frac{1}{\bar L}\bigl(t_ku_k+M_{k+1}\one\bigr),
\qquad M_0=0,
\end{equation}
for some $t_k\in\R$ and $w_k,\theta_k\in\R^r$, where $W_k=\Diag(w_k)$, $\Theta_k=\Diag(\theta_k)$ and $\one\in\R^r$ is the all-ones vector. The $r$ columns of $M_k\in\R^{d\times r}$ are its \emph{momenta}.
\end{definition}

Each momentum is updated with its own weights, and $M_{k+1}$ already contains $u_k$. The matrix $\Theta_0$ plays no role since $M_0=0$. For $r=1$, we write the single momentum as $m_k\in\R^d$, and \eqref{eq:local} reads $m_{k+1}=w_k(u_k+\theta_km_k)$ and $x_{k+1}=x_k-(t_ku_k+m_{k+1})/\bar L$. All methods in this paper have $r=1$ and $t_k=1$ for all $k$. Expanding $M_{k+1}$ gives
\begin{equation}\label{eq:momentum-coeff}
H_{k,i}=\one^\top\Bigl(\prod_{j=i+1}^{k}W_j\Theta_j\Bigr)w_i\quad(i<k),
\qquad
H_{k,k}=t_k+\one^\top w_k.
\end{equation}

\begin{proposition}[momentum form and running sums]\label{prop:momentum}
(i) If $W_k\Theta_k$ is invertible for every $k\ge1$, the momentum form defines the same method as the running-sum form with
\begin{equation}\label{eq:local-additive}
c_k=P_k\one,
\qquad
b_k=P_k^{-1}w_k,
\qquad
d_k=t_k+\one^\top w_k,
\qquad\text{where}\quad
P_k=\prod_{j=1}^{k}W_j\Theta_j,
\end{equation}
and $P_0=I$ is the empty product. \\
(ii) Conversely, if the vectors $b_k$ and $c_k$ have no zero entry, the running-sum form defines the same method as the momentum form with
\[
W_k=\Diag(c_k)\Diag(b_k),
\qquad
\Theta_k=\Diag(c_{k-1})^{-1}\Diag(b_k)^{-1}\quad(k\ge1),
\qquad
t_k=d_k-c_k^\top b_k.
\]
\end{proposition}

\begin{proof}
(i) Since diagonal matrices commute, $P_kP_i^{-1}=\prod_{j=i+1}^{k}W_j\Theta_j$ for $i<k$. Hence \eqref{eq:momentum-coeff} gives $H_{k,i}=\one^\top P_kP_i^{-1}w_i=c_k^\top b_i$ and $H_{k,k}=d_k$, as in \eqref{eq:addfactor}.

(ii) Let $M_{k+1}=S_{k+1}\Diag(c_k)$, the stored vectors rescaled by the read weights. Since $w_k=\Diag(c_k)b_k$ and $\Theta_kW_k=\Diag(c_{k-1})^{-1}\Diag(c_k)$,
\[
M_{k+1}=S_k\Diag(c_k)+u_kb_k^\top\Diag(c_k)=M_k\Theta_kW_k+u_kw_k^\top=\bigl(u_k\one^\top+M_k\Theta_k\bigr)W_k,
\]
and $S_kc_k+d_ku_k=S_{k+1}c_k+(d_k-b_k^\top c_k)u_k=M_{k+1}\one+t_ku_k$.
\end{proof}

The coefficients $c_0$ and $b_{N-1}$ do not affect $H$ and can be chosen freely. For $r=1$, the other coefficients are nonzero when $H$ has no zero entry below the diagonal, since $H_{k,i}=c_kb_i$ for $i<k$.

\subsection{\texorpdfstring{$H$}{H}-duality in momentum form}
\label{sec:duality-momentum}

\begin{proposition}[$H$-dual in momentum form]\label{prop:dual-momentum}
If \eqref{eq:local} is a momentum form of \eqref{eq:Hmap} with $r$ momenta, then the $H$-dual of \eqref{eq:Hmap} has the momentum form
\begin{equation}\label{eq:dual-momentum}
\tilde M_{k+1}=\bigl(u_k\one^\top+\tilde M_k\Theta_{N-k}\bigr)W_{N-1-k},
\qquad
x_{k+1}=x_k-\frac{1}{\bar L}\bigl(t_{N-1-k}u_k+\tilde M_{k+1}\one\bigr),
\end{equation}
with $\tilde M_0=0$, $r$ momenta and $\Theta_N:=0$. The value of $\Theta_N$ plays no role since $\tilde M_0=0$.
\end{proposition}

\begin{proof}
The form \eqref{eq:dual-momentum} has the coefficients $\tilde t_k=t_{N-1-k}$, $\tilde w_k=w_{N-1-k}$ and $\tilde\theta_k=\theta_{N-k}$. Let $\tilde H$ be its coefficient matrix, and let $i<k$. Since diagonal matrices commute, \eqref{eq:momentum-coeff} gives
\[
\tilde H_{k,i}
=\one^\top\Bigl(\prod_{j=i+1}^{k}W_{N-1-j}\Theta_{N-j}\Bigr)W_{N-1-i}\one
=\one^\top\Bigl(\prod_{l=N-1-k}^{N-1-i}W_l\Bigr)\Bigl(\prod_{l=N-k}^{N-1-i}\Theta_l\Bigr)\one.
\]
By \eqref{eq:momentum-coeff} again,
\[
\Hd_{k,i}=H_{N-1-i,\,N-1-k}=\one^\top\Bigl(\prod_{j=N-k}^{N-1-i}W_j\Theta_j\Bigr)W_{N-1-k}\one,
\]
which is the same product. On the diagonal, $\tilde H_{k,k}=t_{N-1-k}+\one^\top w_{N-1-k}=H_{N-1-k,\,N-1-k}=\Hd_{k,k}$. Hence $\tilde H=\Hd$.
\end{proof}

In momentum form, the $H$-dual reverses the sequences $t_k$ and $w_k$, and reverses $\theta_k$ with a one-step shift, because $\theta_k$ acts between steps $k-1$ and $k$. OGM-G \citep{kimfessler2021ogmg}, the $H$-dual of OGM, likewise uses the coefficient sequence of OGM in reverse order \citep{kim2023timereversed}. The $H$-dual coefficients stay bounded when the original ones do, but they may depend on the horizon $N$ even when the original ones do not.

With $r=1$ and $t_k=1$ for every $k$, Proposition~\ref{prop:dual-momentum} shows that palindromic coefficients, $w_k=w_{N-1-k}$ and $\theta_k=\theta_{N-k}$, give a self-dual method, $\Hd=H$. Conversely, $\Hd\neq H$ as soon as $w_0\neq w_{N-1}$, since $H_{0,0}=1+w_0$ by~\eqref{eq:momentum-coeff}, whereas $\Hd_{0,0}=H_{N-1,N-1}=1+w_{N-1}$.

\section{Discovery and minimax optimality of TMM-GD}
\label{sec:tmmgd}

From now on, we consider the minimization of a differentiable function $f:\R^d\to\R$ that is $L$-smooth and $\mu$-strongly convex, with $0<\mu<L$, i.e., for all $x,y\in\R^d$,
\[
f(x)+\ip{\nabla f(x)}{y-x}+\frac\mu2\norm{y-x}^2\le f(y)\le f(x)+\ip{\nabla f(x)}{y-x}+\frac L2\norm{y-x}^2.
\]
We write $\Fcl$ for the class of these functions and $x_\star$ for the unique minimizer of $f$. The vectors $u_k$ of Sections~\ref{sec:memory} and~\ref{sec:recon} are now the gradients $\nabla f(x_k)$, and $\bar L=L$ in \eqref{eq:fsfom}.

\subsection{Decoding the numerical design}
\label{sec:decoding}

\citet{shuwang2026unified} numerically minimize the gradient-to-distance criterion over the coefficient matrices of $5$-step methods. For $\mu/L=1/11$,\footnote{Shu and Wang take $L=1$ and write the strong convexity parameter as $\mu=q_{\text{SW}}/(1+q_{\text{SW}})$, so their value $q_{\text{SW}}=0.1$ corresponds to $\mu/L=1/11$. Their matrix is $6\times6$, with a zero first row for the starting point and a zero last column for the gradient at the output, which is never used. Removing them gives~\eqref{eq:intro-shuwang}.} they report the following coefficient matrix, which we write in the indexing of~\eqref{eq:fsfom} with $\bar L=L$:
\begin{equation}\label{eq:intro-shuwang}
	H\approx
	\begin{pmatrix}
	3.3166&0&0&0&0\\
	0.8684&2.0733&0&0&0\\
	0.3255&0.4024&2.0734&0&0\\
	0.1220&0.1508&0.4024&2.0734&0\\
	0.0987&0.1220&0.3255&0.8684&3.3166
	\end{pmatrix}.
\end{equation} They observe that the matrix is self-antitranspose, i.e., equal to its anti-transpose (Definition~\ref{def:hdual}), consistent with their theory, in which $H$-duality pairs this criterion with itself.

The entries of~\eqref{eq:intro-shuwang} are given to four decimals, so the tails of the columns are known only up to a relative error of order $10^{-4}$. A merge test with a smaller tolerance could therefore reject an exact dependency between the columns because of the rounding, so we run Algorithm~\ref{alg:additive} with $\varepsilon=10^{-3}$, above this level. It returns one stored vector: column~$0$ creates it, and columns $1$, $2$ and $3$ are merged, with residuals at the level of the rounding. Any tolerance between these residuals and $1$ gives the same result. By Theorem~\ref{thm:additive-error}, the resulting coefficient error is below the rounding: each entry of the reconstructed matrix $G$ rounds to the corresponding entry of~\eqref{eq:intro-shuwang}. 

Since $G$ has running-sum memory one and no zero entry below the diagonal, Proposition~\ref{prop:momentum}(ii) writes it in momentum form with one momentum. The coefficients $c_0$ and $b_4$ do not affect $G$, and we choose them so that $t_0=t_4=1$. All the other coefficients are then determined by $G$. Up to the rounding of~\eqref{eq:intro-shuwang}, every step is a gradient step plus momentum ($t_k\approx1$ for every $k$), the weights are palindromic, and $\theta_k$ is constant. They match TMM \citep{vanscoy2018tmm} for $\mu/L=1/11$, written with $\rho=1-\sqrt{\mu/L}$ and $\beta=\rho^2/(2-\rho)$: for $k=1,2,3$, $w_k\approx\rho+\beta$ and $w_k\theta_k\approx\beta$, so these steps run the momentum update $m_{k+1}=\beta m_k+(\rho+\beta)\nabla f(x_k)$ of TMM, and the last weight $w_4\approx\rho/(1-\rho)$ is that of its output rule (Section~\ref{sec:tmmgd-def}). Only the first weight departs from TMM: $w_0\approx\rho/(1-\rho)$ instead of $\rho+\beta$, which sets an initial momentum. The reconstruction procedure and these observations suggest that the numerical design admits a closed form, which we call TMM-GD.

\subsection{TMM-GD}
\label{sec:tmmgd-def}

TMM-GD, stated in Algorithm~\ref{alg:tmmgd}, uses the parameters
\begin{equation}\label{eq:tmm-params}
q:=\frac{\mu}{L},\qquad
\rho:=1-\sqrt q,\qquad
\beta:=\frac{\rho^2}{2-\rho},\qquad
\gamma_\star:=\frac{1+\sqrt q}{2\sqrt q}.
\end{equation}

\begin{algorithm}[H]
\caption{TMM-GD}
\label{alg:tmmgd}
\begin{algorithmic}[1]
\Require starting point $x_0$, horizon $N\ge1$, constants $0<\mu<L$
\Ensure output $z_N$
\State compute $q$, $\rho$, $\beta$ and $\gamma_\star$ from~\eqref{eq:tmm-params}
\State $m_1\gets\dfrac{\rho}{1-\rho}\,\nabla f(x_0)$
\For{$k=1,\dots,N-1$}
\State $x_k\gets x_{k-1}-\dfrac1L\bigl(\nabla f(x_{k-1})+m_k\bigr)$
\State $m_{k+1}\gets\beta m_k+(\rho+\beta)\,\nabla f(x_k)$
\EndFor
\State \Return $z_N\gets x_{N-1}-\dfrac1L\bigl(\nabla f(x_{N-1})+\gamma_\star m_N\bigr)$
\end{algorithmic}
\end{algorithm}

The query points $x_k$ are the points at which the method evaluates gradients. Continuing the loop of Algorithm~\ref{alg:tmmgd} indefinitely gives query points $x_0,x_1,\dots$ and momenta $m_1,m_2,\dots$ that do not depend on $N$. If the method stops after $k\ge1$ gradient evaluations, it returns
\begin{equation}\label{eq:output}
z_k:=x_{k-1}-\frac1L\bigl(\nabla f(x_{k-1})+\gamma_\star m_k\bigr),
\end{equation}
which needs no further gradient. The rule~\eqref{eq:output} does not depend on $N$ either: the horizon only selects which output is returned. In the coefficient matrix of TMM-GD in the form~\eqref{eq:fsfom}, the last row is this output step, so the point denoted $x_N$ there is $z_N$.

With the initial momentum $m_0:=\nabla f(x_0)/(1-\rho)$, line~5 also holds for $k=0$, since $\beta m_0+(\rho+\beta)\nabla f(x_0)=\rho m_0=m_1$. TMM-GD therefore runs a time-invariant recurrence from $x_0$ and $m_0$. This recurrence is TMM: if $\xi_k$ denotes the state of TMM in \citet{vanscoy2018tmm}, the change of variables $x_k=\xi_k+\beta(\xi_k-\xi_{k-1})/(1+\rho)$ and $m_k=-(\rho+\beta)L(\xi_k-\xi_{k-1})/(1+\rho)$ turns the iteration of TMM into lines~4 and~5 of Algorithm~\ref{alg:tmmgd} and its output into~\eqref{eq:output}. TMM-GD is thus TMM started from the state that corresponds to $x_0$ and $m_0$.

The quadratic $f_\mu:=\frac\mu2\norm{\cdot-x_\star}^2$ explains the choice of the initial momentum in Algorithm~\ref{alg:tmmgd}. This quadratic is an extreme function of $\Fcl$, the least curved one: for a prescribed initial gradient norm, it attains the largest initial distance allowed by strong convexity, $\norm{x_0-x_\star}=\norm{\nabla f(x_0)}/\mu$. Extreme functions of this kind are natural worst-cases, so a method optimized for the gradient-to-distance criterion, such as the design~\eqref{eq:intro-shuwang}, tends to improve its behavior on $f_\mu$. On $f_\mu$, the recurrence has a regime in which the momentum is proportional to the error, $m_k=\sqrt{\mu L}\,(x_k-x_\star)$, and both shrink by the factor $\rho$ at every step. Since $\nabla f(x_0)=\mu(x_0-x_\star)$ for $f_\mu$, the initial momentum $m_0=\nabla f(x_0)/(1-\rho)$ equals $\sqrt{\mu L}\,(x_0-x_\star)$: the initial gradient puts the method in this regime from the start, and $x_k-x_\star=\rho^k(x_0-x_\star)$ for every $k$.

The palindromic weights reflect self-duality: for $N\ge2$, Algorithm~\ref{alg:tmmgd} gives the momentum form with $t_k=1$, the constant coefficients $\theta_k=\beta/(\rho+\beta)$, and the palindromic weights $w_0=w_{N-1}=\rho/(1-\rho)$ and $w_k=\rho+\beta$ for $0<k<N-1$, so TMM-GD is its own $H$-dual for every $N$.

For $N=5$ and $\mu/L=1/11$, its coefficient matrix, rounded to four decimals, coincides with~\eqref{eq:intro-shuwang}, up to one rounding error in one entry.

\subsection{Worst-case guarantee of the output}
\label{sec:tmmgd-rate}

\begin{theorem}[rate of the output]\label{thm:rate}
For every $f\in\Fcl$, every $x_0$ and every $N\ge1$, the output $z_N$ of Algorithm~\ref{alg:tmmgd} satisfies
\begin{equation}\label{eq:rate}
\norm{z_N-x_\star}\le\frac{1+\sqrt q}{2\mu}\,\rho^N\,\norm{\nabla f(x_0)},
\end{equation}
with equality for the quadratics $f_\mu$ and $f_L:=\frac L2\norm{\cdot-x_\star}^2$.
\end{theorem}

The bound~\eqref{eq:rate} can be evaluated during a run: $\norm{\nabla f(x_0)}$ is known after the first gradient evaluation, whereas the usual $f(x_0)-f(x_\star)$ and $\norm{x_0-x_\star}$ are not. Since the query points do not depend on $N$, a single run of TMM-GD can stop at the first $k$ with $(1+\sqrt q)\,\rho^k\norm{\nabla f(x_0)}\le2\mu\delta$, for a target accuracy $\delta>0$, and then $\norm{z_k-x_\star}\le\delta$.

\begin{theorem}[interpolation inequalities, \cite{taylorhendrickxglineur2017}]\label{thm:interp}
For points $x_i$ and $x_j$, let
\begin{equation}\label{eq:interp}
\begin{aligned}
I_{ij}:={}&f(x_i)-f(x_j)-\ip{\nabla f(x_j)}{x_i-x_j}-\frac\mu2\norm{x_i-x_j}^2\\
&-\frac{1}{2(L-\mu)}\norm{\nabla f(x_i)-\nabla f(x_j)-\mu(x_i-x_j)}^2.
\end{aligned}
\end{equation}
Every $f\in\Fcl$ satisfies $I_{ij}\ge0$ for all points $x_i, x_j \in \R^d$. Conversely, finitely many points, gradients, and function values that satisfy these inequalities for all pairs are interpolated by a function of $\Fcl$.
\end{theorem}

We call $I_{ij}\ge0$ the \emph{interpolation inequality} between the points $x_i$ and $x_j$. The index $\star$ refers to $x_\star$, where $\nabla f(x_\star)=0$.

The proof is based on a potential that decreases by the constant factor $\rho^2$ at every step:
\begin{equation}\label{eq:potential}
\boxed{\;V_k:=\norm{z_k-x_\star}^2+\frac{L-\mu}{\mu L}\,I_{k-1,\star}\;}\qquad(k\ge1).
\end{equation}
Each decrease $\rho^2V_k-V_{k+1}$ is a nonnegative combination of interpolation inequalities.

\begin{lemma}[potential identities]\label{lem:potential}
The potential~\eqref{eq:potential} satisfies
\begin{align}
V_1&=\Bigl(\frac{L-\mu}{2\mu L}\Bigr)^2\norm{\nabla f(x_0)}^2-\frac{L-\mu}{\mu L}\,I_{\star,0}\label{eq:pot-start}
\intertext{and, for every $k\ge1$,}
\rho^2V_k-V_{k+1}&=\frac{L-\mu}{\mu L}\Bigl[\rho^2I_{k-1,k}+(1-\rho^2)\,I_{\star,k}\Bigr].\label{eq:pot-decay}
\end{align}
\end{lemma}

\begin{proof}
See Appendix~\ref{app:pot}.
\end{proof}

\begin{proof}[Proof of Theorem~\ref{thm:rate}]
By the interpolation inequalities, every $I$ in~\eqref{eq:potential}--\eqref{eq:pot-decay} is nonnegative. Hence $V_1\le\bigl(\frac{L-\mu}{2\mu L}\bigr)^2\norm{\nabla f(x_0)}^2$ by~\eqref{eq:pot-start}, and $V_{k+1}\le\rho^2V_k$ by~\eqref{eq:pot-decay}. Since $V_N\ge\norm{z_N-x_\star}^2$, it follows that $\norm{z_N-x_\star}\le\rho^{N-1}\frac{L-\mu}{2\mu L}\norm{\nabla f(x_0)}$, which is~\eqref{eq:rate} because $(L-\mu)/L=1-q=\rho(1+\sqrt q)$. For $f_\mu$ and $f_L$, these $I$ all vanish, so every inequality above is an equality.
\end{proof}

\begin{remark}[zero initial momentum]\label{rem:initial-momentum}
With zero initial momentum, $m_0=0$, the recurrence starts with too little momentum and reaches the regime of Section~\ref{sec:tmmgd-def} with a larger error. This increases the worst-case distance by the factor $2/(1+\sqrt q)$ at every horizon: the output then satisfies $\norm{z_N-x_\star}\le\rho^N\norm{\nabla f(x_0)}/\mu$, with equality for $f_\mu$. The proof is that of Theorem~\ref{thm:rate}: the potential~\eqref{eq:potential} decreases in the same way, and only its first value changes, to $V_1\le(\rho^2/\mu)\ip{\nabla f(x_0)}{x_0-x_\star}\le\rho^2\norm{\nabla f(x_0)}^2/\mu^2$.  If we instead bound $\ip{\nabla f(x_0)}{x_0-x_\star}$ by $L\norm{x_0-x_\star}^2$, we obtain $\norm{z_N-x_\star}^2\le\rho^{2N}\norm{x_0-x_\star}^2/q$, the criterion and rate of TMM \cite{vanscoy2018tmm}.
\end{remark}

\subsection{Minimax optimality}
\label{sec:tmmgd-minimax}

A method is \emph{linear-span} if each query point $x_k$ lies in $x_0+\operatorname{span}\{\nabla f(x_0),\dots,\nabla f(x_{k-1})\}$ and its output after $N$ gradient evaluations lies in $x_0+\operatorname{span}\{\nabla f(x_0),\dots,\nabla f(x_{N-1})\}$. The coefficients may depend on the observed gradients, so adaptive methods are included. Classical lower bounds \citep{nemirovskyyudin1983,nesterov2004} determine the iteration complexity of smooth strongly convex minimization up to a constant factor. We show that TMM-GD attains the exact minimax value for the gradient-to-distance criterion. Our lower bound builds on the zero-chain construction that \citet{droritaylor2022oracle} use to show that ITEM is optimal for the ratio of final to initial distance. \citet{drori2026weighted} extends this construction and writes its hard functions in closed form, as a quadratic minus a squared distance to a polytope. Our hard function belongs to this family, with parameters chosen for the gradient-to-distance criterion.

To show that no linear-span method has a smaller worst case than TMM-GD, we use a single function $F$ on which every such method has at least this error. Let $d\ge N+1$, and let $e_0,\dots,e_{d-1}$ be the canonical basis of $\R^d$. By Lemma~\ref{lem:hard}, $F\in\mathcal F_{q,1}$, so that $LF\in\Fcl$, and $\nabla F(0)=e_0$, which normalizes the initial gradient. At every point of $\operatorname{span}\{e_0,\dots,e_{k-1}\}$, and not only at particular points, its gradient lies in $\operatorname{span}\{e_0,\dots,e_k\}$, so each gradient evaluation reveals at most one new coordinate, and adaptive choices change the combinations of the observed gradients but cannot leave their span. The minimizer of $F$ has a component along $e_N$, which the output of every linear-span method therefore misses. The parameters below keep $F$ in the class with this chain structure and make the size $\eta_N$ of the missing component as large as possible within this family of functions: this size matches the upper bound of Theorem~\ref{thm:rate}, which makes TMM-GD minimax optimal among all linear-span first-order methods.

Let $q$ and $\rho$ be as in~\eqref{eq:tmm-params}, and set
\begin{align*}
a_0&=1, & \eta_0&=\frac{1+q}{2q},\\
a_i&=\sqrt q\,\eta_i, & \eta_i&=\frac{1-q}{2q^{3/4}}\,\rho^{i-1}\qquad(1\le i<N),\\
a_N&=q\,\eta_N, & \eta_N&=\frac{1-q}{2q}\,\rho^{N-1},
\end{align*}
and $\eta:=\sum_{i=0}^{N}\eta_ie_i$. For $0\le i\le N$, the point $\bar x_i:=-\sum_{t<i}\eta_te_t$ agrees with $-\eta$ on its first $i$ coordinates and vanishes on the others. Let $\mathcal K$ be the convex hull of the points $(a_ie_i-q\bar x_i-\eta)/(1-q)$, $0\le i\le N$, the last of which is $-\eta$. The hard function is
\[
F(x):=\frac12\norm{x+\eta}^2-\frac{1-q}{2}\dist^2(x,\mathcal K),
\]
where $\dist(x,\mathcal K):=\min_{p\in\mathcal K}\norm{x-p}$ is the Euclidean distance from $x$ to $\mathcal K$. Replacing $\mathcal K$ by a single point $p\in\mathcal K$ gives a quadratic of curvature $q$, and $F$ is the maximum of these quadratics over $p\in\mathcal K$.

\begin{lemma}[hard function]\label{lem:hard}
The function $F$ belongs to $\mathcal F_{q,1}$ and its minimizer is $-\eta$. For $0\le i\le N$, $\nabla F(\bar x_i)=a_ie_i$, and in particular, $\nabla F(0)=e_0$. Moreover, for $0\le k<N$, $\nabla F(x)\in\operatorname{span}\{e_0,\dots,e_k\}$ for every $x\in\operatorname{span}\{e_0,\dots,e_{k-1}\}$.
\end{lemma}

\begin{proof}
See Appendix~\ref{app:hard}.
\end{proof}

\begin{theorem}[minimax value]\label{thm:minimax}
Fix $N\ge1$ and $d\ge N+1$. Over all deterministic linear-span methods that use $N$ gradient evaluations, with output $\hat x_N$,
\begin{equation}\label{eq:minimax}
\inf_{\substack{\mathrm{deterministic}\\ \mathrm{linear-span}\\ \mathrm{methods}}}\ \sup_{\substack{f\in\Fcl,\ x_0\\ \nabla f(x_0)\neq0}}
\frac{\norm{\hat x_N-x_\star}^2}{\norm{\nabla f(x_0)}^2}
=\frac{(1+\sqrt q)^2}{4\mu^2}\,\rho^{2N},
\end{equation}
and TMM-GD attains this value.
\end{theorem}

\begin{proof}
TMM-GD is a linear-span method, so Theorem~\ref{thm:rate} gives the upper bound. For the lower bound, take $x_0=0$ and $f=LF$. By Lemma~\ref{lem:hard}, $f\in\Fcl$, $\nabla f(0)=Le_0$, and the minimizer of $f$ is $-\eta$. By the chain property of Lemma~\ref{lem:hard} and induction on $k$, the query points of any linear-span method satisfy $x_k\in\operatorname{span}\{e_0,\dots,e_{k-1}\}$, and its output satisfies $\hat x_N\in\operatorname{span}\{e_0,\dots,e_{N-1}\}$. Since the coordinate of $-\eta$ along $e_N$ is $-\eta_N$,
\[
\frac{\norm{\hat x_N-x_\star}^2}{\norm{\nabla f(0)}^2}\ \ge\ \frac{\eta_N^2}{L^2}=\frac{(1-q)^2}{4q^2L^2}\,\rho^{2N-2}=\frac{(1+\sqrt q)^2}{4\mu^2}\,\rho^{2N}.
\]
\end{proof}

Since $F$ does not depend on the method, the lower bound also holds for randomized methods whose query points and output satisfy the linear-span condition for every realization: on $f=LF$ with $x_0=0$, every realization has a ratio of at least $(1+\sqrt q)^2\rho^{2N}/(4\mu^2)$, and hence so does the expected ratio.

\begin{remark}[methods outside the span]\label{rem:resisting}
For $d\ge2N+1$, the lower bound of Theorem~\ref{thm:minimax} holds for every deterministic first-order method, whose query points and output may leave the span of the observed gradients. The proof is the resisting-oracle argument \citep{nemirovskyyudin1983} in the form used by \citet{droritaylor2022oracle}, with $x_0=0$ and $f=LF(R^\top\cdot)$ for an orthogonal $R$: the oracle chooses the column $Re_k$ when $x_k$ is queried, and $Re_N$ after the output, orthogonal to the earlier columns and to the points queried or returned so far, which imposes at most $2N$ conditions. It applies because $\mathcal K$ and $-\eta$ lie in $E:=\operatorname{span}\{e_0,\dots,e_N\}$, so that $F(y+z)=F(y)+\frac q2\norm{z}^2$ for $y\in E$ and $z\perp E$. With Lemma~\ref{lem:hard}, the answer at $x_k$ then only involves $Re_0,\dots,Re_k$.
\end{remark}

\begin{corollary}[optimal at every horizon]\label{cor:anytime}
A single run of the query recurrence of TMM-GD, which does not depend on $N$, satisfies~\eqref{eq:rate} at every horizon $N\ge1$ through the output $z_N$ of~\eqref{eq:output}. For every horizon $N\le d-1$, this bound is the minimax value~\eqref{eq:minimax}.
\end{corollary}

\begin{proof}
Theorem~\ref{thm:rate} holds for every $N$, and the query points and momenta do not depend on $N$. Theorem~\ref{thm:minimax} applies when $d\ge N+1$.
\end{proof}

Unlike TMM-GD, which runs a time-invariant recurrence with a fixed output rule, exactly optimal methods usually have time-varying coefficients: ITEM only asymptotically corresponds to the stationary TMM \citep{taylordrori2023item}, and the coefficients of ITEM-wf, which is optimal for function values under weighted initial conditions, also vary with the iteration \citep{drori2026weighted}.

\section{\texorpdfstring{$H$}{H}-duality applied to ITEM: ITEM-GG and TMM-GG}
\label{sec:item}

This section applies Proposition~\ref{prop:dual-momentum} to ITEM \citep{taylordrori2023item}, which is optimal for the ratio of final to initial distance \citep{droritaylor2022oracle}. We keep the setting of Section~\ref{sec:tmmgd} and the parameters~\eqref{eq:tmm-params}.

\subsection{ITEM in momentum form}
\label{sec:item-momentum}

ITEM builds its coefficients from the sequence
\begin{equation}\label{eq:item-A}
A_0=0,
\qquad
A_{k+1}=\frac{(1+q)A_k+2\bigl(1+\sqrt{(1+A_k)(1+qA_k)}\bigr)}{(1-q)^2}\quad(k\ge0).
\end{equation}
Its terms are positive for $k\ge1$, and $A_{k+1}\ge A_k/\rho^2$ for every $k\ge0$ \citep{taylordrori2023item}, so the sequence grows at least geometrically.

\begin{proposition}[ITEM in momentum form]\label{prop:item-momentum}
ITEM, with its output denoted $x_N$, runs
\[
m_{k+1}=w_k\bigl(\nabla f(x_k)+\theta_km_k\bigr),
\qquad
x_{k+1}=x_k-\frac1L\bigl(\nabla f(x_k)+m_{k+1}\bigr),
\qquad k=0,\dots,N-1,
\]
from $m_0=0$, with the coefficients
\begin{equation}\label{eq:item-coeffs}
w_k=\frac{2\sqrt{1+A_k}}{(1-q)\sqrt{A_{k+2}}}\ \ (k\le N-2),
\quad
w_{N-1}=\frac{\sqrt{1+A_{N-1}}}{\sqrt{1+qA_N}},
\quad
\theta_k=\frac{A_k}{A_k+(1-q)A_{k+1}}.
\end{equation}
This is the momentum form of Definition~\ref{def:momentum} with one momentum and $t_k=1$ for every $k$.
\end{proposition}

\begin{proof}
See Appendix~\ref{app:item-momentum}.
\end{proof}

Since $A_0=0$, $\theta_0=0$. The last weight $w_{N-1}$ has its own formula because the last step computes the output, as in Algorithm~\ref{alg:tmmgd}. \citet{taylordrori2023item} prove that $\norm{x_N-x_\star}^2\le\norm{x_0-x_\star}^2/(1+qA_N)$ for every $f\in\Fcl$, every $x_0$ and every $N\ge1$, with equality for $f_\mu$ and $f_L$.

\subsection{ITEM-GG}
\label{sec:itemgg}

\citet{shuwang2026unified} pair the ratio $\norm{x_N-x_\star}^2/\norm{x_0-x_\star}^2$ of ITEM with $\norm{\nabla f(x_N)}^2/\norm{\nabla f(x_0)}^2$ under $H$-duality. They implicitly define the $H$-dual of ITEM, which they call ITEM-GG, as the method whose coefficient matrix is the anti-transpose of ITEM's, and prove its rate, but do not state it as an algorithm. By Proposition~\ref{prop:dual-momentum}, ITEM-GG has the momentum form with $t_k=1$, the weights $w_{N-1-k}$ and the coefficients $\theta_{N-k}$, where $w_k$ and $\theta_k$ are the coefficients~\eqref{eq:item-coeffs} of ITEM. Algorithm~\ref{alg:itemgg} states it. Unlike TMM-GD (Section~\ref{sec:tmmgd-def}), ITEM is not self-dual: $w_0<w_{N-1}$ for $N\ge2$ (Appendix~\ref{app:item-momentum}), so ITEM-GG is a different method from ITEM.

\begin{algorithm}[H]
\caption{ITEM-GG}
\label{alg:itemgg}
\begin{algorithmic}[1]
\Require starting point $x_0$, horizon $N\ge1$, constants $0<\mu<L$
\Ensure output $x_N$
\State compute $A_0,\dots,A_N$ from~\eqref{eq:item-A}, and $w_0,\dots,w_{N-1}$ and $\theta_1,\dots,\theta_{N-1}$ from~\eqref{eq:item-coeffs}
\State $m_1\gets w_{N-1}\,\nabla f(x_0)$
\For{$k=1,\dots,N-1$}
\State $x_k\gets x_{k-1}-\dfrac1L\bigl(\nabla f(x_{k-1})+m_k\bigr)$
\State $m_{k+1}\gets w_{N-1-k}\bigl(\nabla f(x_k)+\theta_{N-k}\,m_k\bigr)$
\EndFor
\State \Return $x_N\gets x_{N-1}-\dfrac1L\bigl(\nabla f(x_{N-1})+m_N\bigr)$
\end{algorithmic}
\end{algorithm}

\begin{theorem}[rate of ITEM-GG, \cite{shuwang2026unified}]\label{thm:itemgg}
For every $f\in\Fcl$, every $x_0$ and every $N\ge1$, the output $x_N$ of Algorithm~\ref{alg:itemgg} satisfies
\[
\norm{\nabla f(x_N)}^2\le\frac{\norm{\nabla f(x_0)}^2}{1+qA_N},
\]
with equality for $f_\mu$ and $f_L$.
\end{theorem}

Instead of transferring a Performance Estimation certificate of ITEM through $H$-duality, we prove Theorem~\ref{thm:itemgg} directly with a potential in the iterates and momenta of ITEM-GG, as \citet{taylordrori2023item} do for ITEM. Each of its decreases is a nonnegative combination of interpolation inequalities.

The potential~\eqref{eq:potential} of TMM-GD is written against the minimizer $x_\star$, but the bound of Theorem~\ref{thm:itemgg} involves no minimizer. For a point $x$, let $x^+:=x-\nabla f(x)/L$. It minimizes the quadratic upper bound $f(x)+\ip{\nabla f(x)}{\cdot-x}+\frac L2\norm{\cdot-x}^2$ of $f$, where this bound has the gradient $0$ and the value $f(x)-\norm{\nabla f(x)}^2/(2L)$. With this gradient and this value, we call $x^+$ the \emph{virtual minimizer} of $x$. For smooth convex functions, analyses of OGM-G \citep{diakonikolaswang2022,leeparkryu2021} and the $H$-dual energies of \citet{kim2023timereversed} are anchored at the last iterate, and these analyses of OGM-G also involve the gradient step of this iterate, with which \citet{diakonikolaswang2022} rewrite the interpolation inequalities anchored at it. Expanding~\eqref{eq:interp} gives $I_{x^+,x_j}=I_{x,x_j}$ for every point $x_j$ (Appendix~\ref{app:identities}), also when $\mu=0$, so the anchor becomes its virtual minimizer, a point with zero gradient against which certificates written for $x_\star$ can be reused (Section~\ref{sec:tmmgg}).

The potential of ITEM-GG replaces $x_\star$ by the virtual minimizer $x_N^+$ of the output: in place of the interpolation inequalities with $x_\star$ as the first point, its certificate uses the inequalities $I_{x_N^+,x_k}=I_{N,k}$. This choice fits the criterion: since $\norm{\nabla f(x_N)}=L\norm{x_N-x_N^+}$, the criterion measures the distance from the output to its virtual minimizer, as that of TMM-GD measures the distance from the output to $x_\star$. The identities of the potential are local in the iteration index, but, through $x_N^+$, the certificate depends on $x_N$. The potential is, for $1\le k\le N$,
\begin{equation}\label{eq:itemgg-potential}
\boxed{\;\begin{aligned}\Phi_k:={}&\norm{\nabla f(x_k)+\theta_{N-k}m_k}^2+L\,\ip{\nabla f(x_k)-(1-2\theta_{N-k})\,m_k}{x_k^+-x_N^+}\\&-(L-\mu)\,I_{N,k}.\end{aligned}\;}
\end{equation}
At $k=N$, the inner product and $I_{N,N}$ vanish and $\theta_0=0$, so $\Phi_N=\norm{\nabla f(x_N)}^2$.

\begin{lemma}[potential identities of ITEM-GG]\label{lem:itemgg-potential}
The potential~\eqref{eq:itemgg-potential} satisfies
\begin{align}
\Phi_1&=\frac{(1-q)^2A_N}{4(1+qA_N)}\,\norm{\nabla f(x_0)}^2-(L-\mu)\bigl(I_{0,1}+I_{N,0}\bigr)\label{eq:itemgg-start}
\intertext{and, for $1\le k\le N-1$,}
\frac{A_{N-k}}{A_{N-k+1}}\,\Phi_k-\Phi_{k+1}&=(L-\mu)\Bigl[I_{k,k+1}+\Bigl(1-\frac{A_{N-k}}{A_{N-k+1}}\Bigr)I_{N,k}\Bigr].\label{eq:itemgg-decay}
\end{align}
\end{lemma}

\begin{proof}
See Appendix~\ref{app:itemgg-potential}.
\end{proof}

\begin{proof}[Proof of Theorem~\ref{thm:itemgg}]
By the interpolation inequalities, every $I$ in~\eqref{eq:itemgg-start} and~\eqref{eq:itemgg-decay} is nonnegative. Since $0<A_{N-k}<A_{N-k+1}$, it follows that $\Phi_{k+1}\le\frac{A_{N-k}}{A_{N-k+1}}\Phi_k$ for $1\le k\le N-1$. The product of these factors telescopes to $A_1/A_N$, and $A_1=4/(1-q)^2$. Since $\Phi_N=\norm{\nabla f(x_N)}^2$, \eqref{eq:itemgg-start} gives
\[
\norm{\nabla f(x_N)}^2\le\frac{A_1}{A_N}\,\Phi_1\le\frac{A_1}{A_N}\,\frac{(1-q)^2A_N}{4(1+qA_N)}\,\norm{\nabla f(x_0)}^2=\frac{\norm{\nabla f(x_0)}^2}{1+qA_N}.
\]
For $f_\mu$ and $f_L$, these $I$ all vanish, so every inequality above is an equality.
\end{proof}

For $N=1$, ITEM-GG is the gradient step with step size $2/(L+\mu)$, which is optimal for this ratio. We show that, for every $N$, in dimension $d\ge N+1$, no nonadaptive linear-span method does better, and that ITEM-GG is the only one that attains the bound.

\begin{theorem}[minimax optimality of ITEM-GG]\label{thm:itemgg-minimax}
Fix $N\ge1$ and $d\ge N+1$. For every nonadaptive linear-span method~\eqref{eq:fsfom},
\[
\sup_{\substack{f\in\Fcl,\ x_0\\ \nabla f(x_0)\neq0}}\frac{\norm{\nabla f(x_N)}^2}{\norm{\nabla f(x_0)}^2}\ \ge\ \frac1{1+qA_N},
\]
with equality only for ITEM-GG. Hence ITEM-GG is the unique minimax optimal nonadaptive linear-span method for this ratio.
\end{theorem}

\begin{proof}
See Appendix~\ref{app:itemgg-minimax}.
\end{proof}

\begin{corollary}[dependence on the horizon]\label{cor:horizon}
Let $2\le N_1<N_2$ and $d\ge N_2+1$. No single run of a nonadaptive linear-span method attains the minimax value of Theorem~\ref{thm:itemgg-minimax} at both horizons $N_1$ and $N_2$, even if its output at $N_1$ differs from its query point $x_{N_1}$.
\end{corollary}

\begin{proof}
By Theorem~\ref{thm:itemgg-minimax}, attaining the minimax value at a horizon $N\ge2$ requires the first query point $x_1=x_0-\sqrt{A_N/(1+qA_N)}\,\nabla f(x_0)/L$ of ITEM-GG (Algorithm~\ref{alg:itemgg} and Appendix~\ref{app:item-momentum}), which differs for $N_1$ and $N_2$, since $A_N/(1+qA_N)$ increases strictly with $N$.
\end{proof}

\citet{shuwang2026unified} describe the rate of Theorem~\ref{thm:itemgg} as the conjectured optimal gradient-to-gradient guarantee, and Theorem~\ref{thm:itemgg-minimax} proves this optimality among nonadaptive linear-span methods. ITEM and ITEM-GG thus relate as OGM and its $H$-dual OGM-G do for smooth convex functions: OGM is optimal for the ratio $(f(x_N)-f(x_\star))/\norm{x_0-x_\star}^2$ \citep{drori2017exact}, and OGM-G is the unique optimal nonadaptive linear-span method for the dual ratio $\norm{\nabla f(x_N)}^2/(f(x_0)-f(x_\star))$ \citep{grimmerjopark2026}. In both pairs, the first method depends on the horizon $N$ only through its last step, which computes the output, whereas every step of its $H$-dual depends on $N$, including the first one, since $H$-duality reverses time. In other settings, $H$-duality does not always preserve optimality: for nonexpansive fixed-point problems, \citet{yoongrimmer2026} determine when it does. For these problems, optimal methods are not unique \citep{yoonkimsuhryu2024,yoonryugrimmer2025}, whereas Theorem~\ref{thm:itemgg-minimax} shows that ITEM-GG is the only optimal method.

The proof of Theorem~\ref{thm:itemgg-minimax} adapts the strategy of \citet{grimmerjopark2026} for OGM-G to $\Fcl$ and the gradient-to-gradient criterion. By Performance Estimation duality, the worst case of a method is the infimum of the values $\tau$ certified by nonnegative combinations of interpolation inequalities between its iterates. Eliminating the rows of the method in reverse order, each together with a Schur complement of the multipliers, reduces these certificates to scalars $\delta_0,\dots,\delta_{N-1}$, terminal weights attached to $x_N$. Appendix~\ref{app:itemgg-minimax} defines these scalars and the increasing map $\psi$ below. Every certificate with positive multipliers satisfies
\[
\begin{gathered}
\delta_{N-1}\ge\frac1{A_1},
\qquad
\delta_{m-1}\ge\psi(\delta_m)\quad(1\le m\le N-1),
\qquad
\psi\Bigl(\frac1{A_j}\Bigr)=\frac1{A_{j+1}}\quad(j\ge1),\\
\tau\ge\frac{\delta_0}{q+\delta_0}\ge\frac1{1+qA_N},
\end{gathered}
\]
where the last inequality uses $\delta_0\ge1/A_N$, which follows from the first line by induction. The identity $\psi(1/A_j)=1/A_{j+1}$ is the recursion~\eqref{eq:item-A} written for $1/A_j$: in this adaptation, the scalar reduction forces the ITEM sequence instead of that of OGM, and the equality cases of the chain determine ITEM-GG. As in \citet{grimmerjopark2026}, uniqueness follows from these equality cases through limits of certificates. Since our criterion does not involve $x_\star$, $d\ge N+1$ suffices instead of their $d\ge N+2$. Unlike Theorem~\ref{thm:minimax}, Theorem~\ref{thm:itemgg-minimax} only concerns nonadaptive methods: its lower bound depends on the method through $H$, whereas the hard function of Theorem~\ref{thm:minimax} serves every linear-span method. The potentials show a similar contrast: in that of TMM-GD, the minimizer $x_\star$ only depends on $f$, whereas in that of ITEM-GG, its role is played by the virtual minimizer $x_N^+$, which also depends on the method.

\subsection{TMM-GG}
\label{sec:tmmgg}

ITEM-GG must fix the horizon before it computes $x_1$. The guarantee of ITEM-GG only concerns $x_N$, so stopping a run earlier, or continuing it, is not covered. This makes ITEM-GG impractical when the number of iterations is not known in advance. The obstacle is this dependence on $N$, not the variation of the coefficients with $k$, which ITEM also has. By Corollary~\ref{cor:horizon}, no nonadaptive linear-span method avoids it while remaining minimax optimal at two horizons $N_1,N_2\ge2$.

We therefore remove the dependence on $N$ by freezing the coefficients of ITEM-GG at their limits as $N\to\infty$, and we show that the resulting method, TMM-GG, keeps the accelerated rate at every horizon of a single run, at the price of a factor between $1/(4q)$ and $1/q$ in the squared bound with respect to ITEM-GG (Remark~\ref{rem:price}).

Since $H$-duality reverses time, the first steps of ITEM-GG use the last coefficients of ITEM, and ITEM corresponds asymptotically to TMM \citep{taylordrori2023item}: $A_k/A_{k+1}\to\rho^2$ as $k\to\infty$. Inserting this limit into~\eqref{eq:item-coeffs}, with $1-q=\rho(2-\rho)$ and $\rho+\beta=2\rho/(2-\rho)$, gives $w_{N-1}\to\rho/(1-\rho)$, $w_{N-1-k}\to\rho+\beta$ and $\theta_{N-k}\to\rho/2$ as $N\to\infty$, for every fixed $k\ge1$. The method with these limits as coefficients, which we call TMM-GG, does not depend on $N$. \citet{upadhyaya2026proxitem} also study the stationary limit of an optimal method for composite problems.

In momentum form, TMM-GG has $t_k=1$, the weight $\rho/(1-\rho)$ at step~$0$, and the weight $\rho+\beta$ and the coefficient $\beta/(\rho+\beta)=\rho/2$ at steps $1,\dots,N-1$. These are the coefficients of TMM-GD, except for the last step: TMM-GG runs the loop of Algorithm~\ref{alg:tmmgd} and returns the query point $x_N$ in place of $z_N$, so its last step uses $m_N$ instead of $\gamma_\star m_N$. A single run of this loop thus serves both criteria at every horizon: $z_N$ for the distance to $x_\star$ (Corollary~\ref{cor:anytime}) and $x_N$ for the gradient norm (Theorem~\ref{thm:tmmgg} below). By Proposition~\ref{prop:dual-momentum}, TMM-GG is also the $H$-dual of TMM with zero initial momentum: reversing the weights of TMM-GG puts the weight $\rho+\beta$, which corresponds to $m_0=0$, at step~$0$, and the weight $\rho/(1-\rho)=\gamma_\star(\rho+\beta)$, which returns $z_N$, at the last step.

\begin{theorem}[rate of TMM-GG]\label{thm:tmmgg}
For every $f\in\Fcl$, every $x_0$ and every $N\ge1$, the query point $x_N$ of Algorithm~\ref{alg:tmmgd} satisfies
\[
\norm{\nabla f(x_N)}\le\frac{\rho^N}{\sqrt q}\,\norm{\nabla f(x_0)},
\]
with equality for $f_L$.
\end{theorem}

\begin{proof}
See Appendix~\ref{app:momentum-bound}.
\end{proof}

Since $\nabla f(x_N)$ is a combination of $m_N$ and $m_{N+1}$, the proof bounds the momenta. To bound $m_k$, $k\ge1$, it applies the potential~\eqref{eq:potential} to the loop that TMM-GG shares with TMM-GD, with $x_\star$ replaced by the virtual minimizer $x_{k-1}^+$ of Section~\ref{sec:itemgg}. The points $z_k$ of~\eqref{eq:output} in this potential only serve the analysis, since TMM-GG returns $x_N$. The proof of Theorem~\ref{thm:rate} still works with $x_{k-1}^+$ in place of $x_\star$: the identities of Lemma~\ref{lem:potential} only use the zero gradient at $x_\star$, the interpolation inequalities with $x_\star$ as first point become $I_{x_{k-1}^+,x_j}=I_{x_{k-1},x_j}\ge0$, and the one with $x_\star$ as second point becomes $I_{x_{k-1},x_{k-1}^+}=0$ (Appendix~\ref{app:momentum-bound}). It therefore bounds the distance from $z_k$ to $x_{k-1}^+$, which is $\gamma_\star\norm{m_k}/L$ since $z_k=x_{k-1}^+-\gamma_\star m_k/L$ by~\eqref{eq:output}. Hence $\norm{m_k}\le\rho^k\norm{\nabla f(x_0)}/\sqrt q$ for every $k\ge1$, with equality for $f_\mu$ and $f_L$. The potential of TMM-GD thus certifies both methods, written against the minimizer for the distance and against virtual minimizers for the gradient norm.

Theorem~\ref{thm:tmmgg} is an anytime guarantee in the sense of \citet{kornowski2024open} and \citet{zhang2025anytime}, with the accelerated rate $\rho^{2N}$. Strong convexity is essential here: for smooth convex functions, \citet{vernimmen2026anytime} show that no nonadaptive linear-span method whose coefficients do not depend on the horizon reduces $\norm{\nabla f(x_N)}^2/(f(x_0)-f(x_\star))$ at an accelerated rate for every horizon, a weak form of a conjecture of \citet{diakonikolaswang2022}, so the horizon dependence of OGM-G cannot be removed without losing the accelerated rate of its last iterate at infinitely many horizons.

\begin{figure}[h]
\centering
\definecolor{priceA}{HTML}{2A78D6}
\definecolor{priceB}{HTML}{EB6834}
\definecolor{priceC}{HTML}{1BAF7A}
\begin{tikzpicture}
\begin{semilogyaxis}[
  width=0.82\textwidth,height=0.42\textwidth,
  xmin=0,xmax=41,ymin=1,ymax=120,
  xtick={1,5,10,15,20,25,30,35,40},
  ytick={1,2,5,10,20,50,100},yticklabels={$1$,$2$,$5$,$10$,$20$,$50$,$100$},
  xlabel={$N$},ylabel={ratio of the squared bounds},
  legend style={at={(0.5,1.03)},anchor=south,legend columns=5,font=\small,/tikz/every even column/.append style={column sep=0.8em}},legend cell align=left,
  tick label style={font=\small},label style={font=\small},
]
\addlegendimage{priceA,thick}\addlegendentry{$q=1/4$}
\addlegendimage{priceB,thick}\addlegendentry{$q=1/11$}
\addlegendimage{priceC,thick}\addlegendentry{$q=1/100$}
\addlegendimage{black,thick,dashed}\addlegendentry{output $x_N$}
\addlegendimage{black,thick,densely dotted}\addlegendentry{output $z^G_N$}
\addplot[priceA,thick,dashed,forget plot] coordinates {(1,2.7778) (2,2.5147) (3,2.4514) (4,2.4357) (5,2.4318) (6,2.4308) (7,2.4306) (8,2.4305) (9,2.4305) (10,2.4305) (11,2.4305) (12,2.4305) (13,2.4305) (14,2.4305) (15,2.4305) (16,2.4305) (17,2.4305) (18,2.4305) (19,2.4305) (20,2.4305) (21,2.4305) (22,2.4305) (23,2.4305) (24,2.4305) (25,2.4305) (26,2.4305) (27,2.4305) (28,2.4305) (29,2.4305) (30,2.4305) (31,2.4305) (32,2.4305) (33,2.4305) (34,2.4305) (35,2.4305) (36,2.4305) (37,2.4305) (38,2.4305) (39,2.4305) (40,2.4305)};
\addplot[priceA,thick,densely dotted,forget plot] coordinates {(1,1.5625) (2,1.4145) (3,1.3789) (4,1.3701) (5,1.3679) (6,1.3673) (7,1.3672) (8,1.3672) (9,1.3672) (10,1.3672) (11,1.3672) (12,1.3672) (13,1.3672) (14,1.3672) (15,1.3672) (16,1.3672) (17,1.3672) (18,1.3672) (19,1.3672) (20,1.3672) (21,1.3672) (22,1.3672) (23,1.3672) (24,1.3672) (25,1.3672) (26,1.3672) (27,1.3672) (28,1.3672) (29,1.3672) (30,1.3672) (31,1.3672) (32,1.3672) (33,1.3672) (34,1.3672) (35,1.3672) (36,1.3672) (37,1.3672) (38,1.3672) (39,1.3672) (40,1.3672)};
\addplot[priceB,thick,dashed,forget plot] coordinates {(1,7.7281) (2,6.3856) (3,5.785) (4,5.5044) (5,5.3704) (6,5.3057) (7,5.2743) (8,5.259) (9,5.2516) (10,5.248) (11,5.2462) (12,5.2453) (13,5.2449) (14,5.2447) (15,5.2446) (16,5.2446) (17,5.2445) (18,5.2445) (19,5.2445) (20,5.2445) (21,5.2445) (22,5.2445) (23,5.2445) (24,5.2445) (25,5.2445) (26,5.2445) (27,5.2445) (28,5.2445) (29,5.2445) (30,5.2445) (31,5.2445) (32,5.2445) (33,5.2445) (34,5.2445) (35,5.2445) (36,5.2445) (37,5.2445) (38,5.2445) (39,5.2445) (40,5.2445)};
\addplot[priceB,thick,densely dotted,forget plot] coordinates {(1,3.2727) (2,2.7042) (3,2.4499) (4,2.331) (5,2.2743) (6,2.2469) (7,2.2336) (8,2.2271) (9,2.224) (10,2.2224) (11,2.2217) (12,2.2213) (13,2.2211) (14,2.221) (15,2.221) (16,2.221) (17,2.221) (18,2.221) (19,2.221) (20,2.221) (21,2.221) (22,2.221) (23,2.221) (24,2.221) (25,2.221) (26,2.221) (27,2.221) (28,2.221) (29,2.221) (30,2.221) (31,2.221) (32,2.221) (33,2.221) (34,2.221) (35,2.221) (36,2.221) (37,2.221) (38,2.221) (39,2.221) (40,2.221)};
\addplot[priceC,thick,dashed,forget plot] coordinates {(1,84.306) (2,72.787) (3,64.166) (4,57.626) (5,52.61) (6,48.728) (7,45.701) (8,43.325) (9,41.45) (10,39.963) (11,38.78) (12,37.835) (13,37.079) (14,36.473) (15,35.985) (16,35.593) (17,35.277) (18,35.022) (19,34.816) (20,34.65) (21,34.516) (22,34.407) (23,34.319) (24,34.248) (25,34.191) (26,34.144) (27,34.107) (28,34.076) (29,34.051) (30,34.031) (31,34.015) (32,34.002) (33,33.991) (34,33.983) (35,33.976) (36,33.97) (37,33.966) (38,33.962) (39,33.959) (40,33.957)};
\addplot[priceC,thick,densely dotted,forget plot] coordinates {(1,25.503) (2,22.018) (3,19.41) (4,17.432) (5,15.914) (6,14.74) (7,13.824) (8,13.106) (9,12.538) (10,12.089) (11,11.731) (12,11.445) (13,11.216) (14,11.033) (15,10.886) (16,10.767) (17,10.671) (18,10.594) (19,10.532) (20,10.482) (21,10.441) (22,10.408) (23,10.382) (24,10.36) (25,10.343) (26,10.329) (27,10.317) (28,10.308) (29,10.301) (30,10.294) (31,10.29) (32,10.286) (33,10.282) (34,10.28) (35,10.278) (36,10.276) (37,10.275) (38,10.273) (39,10.273) (40,10.272)};
\end{semilogyaxis}
\end{tikzpicture}
\caption{Factor by which the squared bound of TMM-GG exceeds that of ITEM-GG, with the output $x_N$ (dashed, Theorem~\ref{thm:tmmgg}) and with the output $z^G_N$ (dotted, Theorem~\ref{thm:tmmgg-out}).}
\label{fig:price}
\end{figure}
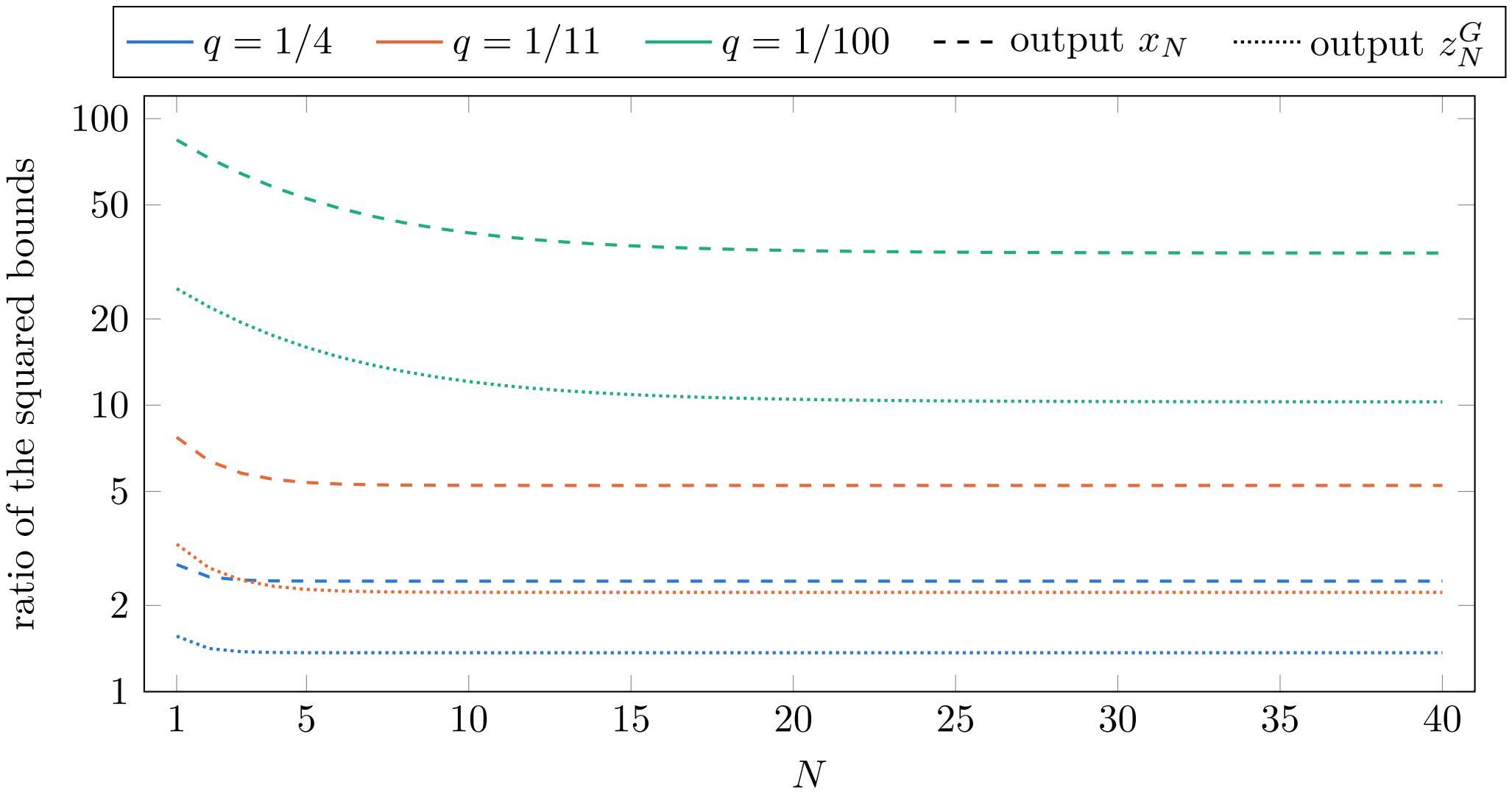

The certificate of TMM-GG also suggests that the method can be improved: the potential of TMM-GD can be written against the virtual minimizer of any point, not only of a query point (Appendix~\ref{app:tmmgg-out}), so it certifies other outputs as well, and we use this freedom to change the output. Theorem~\ref{thm:tmmgg} is an equality for $f_L$, but for $f_\mu$, $x_N-x_\star=\rho^N(x_0-x_\star)$ (Section~\ref{sec:tmmgd-def}), so $\norm{\nabla f(x_N)}=\rho^N\norm{\nabla f(x_0)}$ is smaller by the factor $\sqrt q$. Shortening the last momentum step balances the two quadratics: the point
\begin{equation}\label{eq:output-gg}
z^G_N:=x_{N-1}-\frac1L\bigl(\nabla f(x_{N-1})+\gamma_Gm_N\bigr),
\qquad
\gamma_G:=\frac{\sqrt q\,(1+\sqrt q)}{1+q},
\end{equation}
satisfies $\norm{\nabla f(z^G_N)}=\frac{1+\sqrt q}{1+q}\,\rho^N\norm{\nabla f(x_0)}$ for both $f_\mu$ and $f_L$.

\begin{theorem}[rate of the output $z^G_N$]\label{thm:tmmgg-out}
For every $f\in\Fcl$, every $x_0$ and every $N\ge1$, the point $z^G_N$ of~\eqref{eq:output-gg}, formed from the iterates and momenta of Algorithm~\ref{alg:tmmgd}, satisfies
\[
\norm{\nabla f(z^G_N)}^2\le\frac{(1+\sqrt q)^2}{4q}\,\rho^{2N}\,\norm{\nabla f(x_0)}^2.
\]
\end{theorem}

\begin{proof}
See Appendix~\ref{app:tmmgg-out}.
\end{proof}

The bound of Theorem~\ref{thm:tmmgg-out} is smaller than that of Theorem~\ref{thm:tmmgg} by the factor $(1+\sqrt q)^2/4$. For $2\le N\le6$ and the values of $q$ in Figure~\ref{fig:price}, Performance Estimation, solved with PEPit \citep{goujaud2024pepit}, gives a worst case within $5\%$ of this bound. Since $\gamma_G$ does not depend on $N$, Theorem~\ref{thm:tmmgg-out} holds at every horizon of a single run, and, like $z_N$, the point $z^G_N$ needs no further gradient evaluation. The three outputs of the loop of Algorithm~\ref{alg:tmmgd} considered so far are points of the family
\[
z_N(\gamma):=x_{N-1}^+-\frac\gamma L\,m_N,
\]
namely $z_N(1)=x_N$, $z_N(\gamma_\star)=z_N$ and $z_N(\gamma_G)=z^G_N$. The common query recurrence thus serves both criteria, which only enter through the coefficient of the last momentum step.

\begin{remark}[cost of horizon independence]\label{rem:price}
The squared bound of TMM-GG in Theorem~\ref{thm:tmmgg} exceeds that of ITEM-GG in Theorem~\ref{thm:itemgg} by the factor $R_N:=\rho^{2N}(1+qA_N)/q$. This factor decreases with $N$ and satisfies $1/(4q)<R_N<1/q$ for every $N$ (Appendix~\ref{app:item-identities}). Since the bound of TMM-GG decreases by the factor $\rho^2<e^{-2\sqrt q}$ at each iteration, TMM-GG guarantees, after $N+\lceil\ln(1/q)/(2\sqrt q)\rceil$ iterations, at least the reduction of the gradient norm that ITEM-GG guarantees after $N$ iterations. Figure~\ref{fig:price} shows $R_N$ and the factor $R_N(1+\sqrt q)^2/4$ of the output $z^G_N$.
\end{remark}

\FloatBarrier

\backmatter

\section*{Statements and Declarations}

\noindent\textbf{Funding.} Pierre Vernimmen is funded by the Fonds Sp\'eciaux de Recherche (FSR) of UCLouvain.
A. Taylor is supported by the European Union (ERC grant
CASPER 101162889).
B. Goujaud is supported by a Hi!Paris
chair. The French government also partly funded this work under the management of the Agence
Nationale de la Recherche as part of the France 2030 program, references ANR-23-IACL-0008
“PR[AI]RIE-PSAI”, ANR-23-PEIA-005 (REDEEM project), and ANR-23-IACL-0005.

\noindent\textbf{Competing interests.} The authors have no competing interests to declare.

\noindent\textbf{Use of generative AI.} The authors used LLMs (Claude Opus 5.5 and ChatGPT Astra 6) during the preparation of this work to draft, revise and check parts of the text, of the derivations, and of the proofs. They also used them to search the literature and check references, to write code for the numerical checks, and to prepare the LaTeX source. The authors verified the output of these tools and take full responsibility for the content of this article.

\appendix

\section{Memory and running-sum memory}
\label{app:general}

\begin{proof}[Proof of Theorem~\ref{thm:addmemory}(i) and (iii)]
Let $j\ge k>i$. By \eqref{eq:general-factor}, $H_{j,i}$ is the product of the row vector $b_i^\top A_{i+1}\cdots A_{k-1}$, which depends only on $i$, and the column vector $A_k\cdots A_{j-1}c_j$, which depends only on $j$. Hence $K_k(H)$ has rank at most $r$, and $r(H)\ge\max_k\rank K_k(H)$.

Conversely, let $\bar r$ be this maximum. For $k=1,\dots,N-1$, let the rows of $R_k\in\R^{\bar r\times k}$ be a basis of the row space of $K_k(H)$, completed by zero rows, and write $K_k(H)=O_kR_k$. Set $S_k=[\,u_0\ \cdots\ u_{k-1}\,]R_k^\top$, and let $c_k^\top$ be the first row of $O_k$. Then $S_kc_k=\sum_{i<k}H_{k,i}u_i$, so step $k$ of \eqref{eq:Hmap} reads $x_{k+1}=x_k-(S_kc_k+H_{k,k}u_k)/\bar L$. The first $k$ columns of $K_{k+1}(H)$ form a submatrix of $K_k(H)$, namely its rows $k+1,\dots,N-1$. Restricted to these columns, the rows of $R_{k+1}$ therefore lie in the row space of $R_k$, so $R_{k+1}=[\,A_k^\top R_k\ \ b_k\,]$ for some $A_k\in\R^{\bar r\times\bar r}$ and $b_k\in\R^{\bar r}$. Then $S_{k+1}=S_kA_k+u_kb_k^\top$. With $b_0=R_1$ and $d_k=H_{k,k}$, this gives \eqref{eq:general} with $\bar r$ stored vectors. Finally, running sums are general recurrences, so $r(H)\le\rsum(H)$.

For the equality case, deleting the first row of $K_k(H)$ gives $K_{k+1}^-(H)$, and appending the tail of column $k$ then gives $K_{k+1}(H)$. The first operation lowers the rank by one exactly when step $k$ forgets, and the second raises it by one exactly when column $k$ is new. Since $\rank K_1(H)$ is one if column $0$ is new and zero otherwise, $\rank K_m(H)$ is the number of new columns $i<m$ minus the number of steps $k<m$ that forget. By Theorem~\ref{thm:addmemory}(ii), $\rsum(H)$ is the number of all new columns. Hence $\rank K_m(H)=\rsum(H)$ if and only if no step $k<m$ forgets and no column $i\ge m$ is new. Such an $m$ exists if and only if $i<k$ whenever column $i$ is new and step $k$ forgets.
\end{proof}

\begin{corollary}[generic values]\label{cor:generic}
For almost every lower-triangular $H\in\R^{N\times N}$, $r(H)=\rsum(H)=\lfloor N/2\rfloor$.
\end{corollary}

\begin{proof}
The entries of $K_k(H)$ and $K_k^-(H)$ lie strictly below the diagonal, so, for almost every $H$, these blocks have full rank: $\rank K_k(H)=\min(k,N-k)$ and $\rank K_k^-(H)=\min(k-1,N-k)$. The term of~\eqref{eq:radd} with index $k$ is therefore $1$ if $k\le N/2$ and $0$ otherwise, so $\rsum(H)=\lfloor N/2\rfloor$, and Theorem~\ref{thm:addmemory}(i) gives $r(H)=\max_k\min(k,N-k)=\lfloor N/2\rfloor$.
\end{proof}

\section{Proofs for Section~\ref{sec:tmmgd}}
\label{app:tmmgd}

\subsection{Interpolation identities}
\label{app:identities}

The proofs below use the following identities between the quantities~\eqref{eq:interp}, for $L=1$ and $\mu=q$. Two combinations of $x_i$ and $\nabla f(x_i)$ appear in all the computations:
\[
g^\mu_i:=\nabla f(x_i)-q(x_i-x_\star),
\qquad
g^L_i:=(x_i-x_\star)-\nabla f(x_i),
\]
with $g^\mu_\star=g^L_\star=0$. With $\phi_i:=f(x_i)-f(x_\star)-\frac q2\norm{x_i-x_\star}^2$, the definition~\eqref{eq:interp} reads
\[
I_{ij}=\phi_i-\phi_j-\ip{g^\mu_j}{x_i-x_j}-\frac{1}{2(1-q)}\norm{g^\mu_i-g^\mu_j}^2,
\]
and $\phi_\star=0$. For all indices $i$ and $j$, including $\star$, it follows that
\begin{equation}\label{eq:twopoint}
(1-q)\bigl(I_{i\star}+I_{\star j}-I_{ij}\bigr)=\ip{g^\mu_j}{(1-q)(x_i-x_\star)-g^\mu_i}=\ip{g^\mu_j}{g^L_i}.
\end{equation}
Applying~\eqref{eq:twopoint} to the pairs $(i,j)=(a,a)$, $(a,b)$, $(c,a)$ and $(c,b)$, adding the first and last identities and subtracting the other two gives, for all indices $a$, $b$ and $c$,
\begin{equation}\label{eq:threepoint}
(1-q)\bigl(I_{ab}+I_{ca}-I_{cb}\bigr)=\ip{g^\mu_a-g^\mu_b}{g^L_a-g^L_c}.
\end{equation}
In both identities, the indices may refer to any points of $\R^d$, not only to iterates.

For a point $x$, the virtual minimizer of Section~\ref{sec:itemgg} is the gradient step $x^+=x-\nabla f(x)$, with the gradient $0$ and the value $f(x)-\frac12\norm{\nabla f(x)}^2$. Expanding~\eqref{eq:interp} gives, for every point $x_j$,
\begin{equation}\label{eq:virtual}
I_{x^+,x_j}=I_{x,x_j},
\qquad
I_{x_j,x^+}=I_{x_j,x}+\ip{\nabla f(x)}{x_j^+-x^+},
\end{equation}
and in particular $I_{x,x^+}=0$.

\subsection{Potential identities}
\label{app:pot}
 
\begin{proof}[Proof of Lemma~\ref{lem:potential}]
Replacing $f$ by $f/L$ leaves the iterates of Algorithm~\ref{alg:tmmgd} and both sides of~\eqref{eq:pot-start} and~\eqref{eq:pot-decay} unchanged, and turns $(\mu,L)$ into $(q,1)$. We may therefore assume $L=1$ and $\mu=q$, and use the notation of Appendix~\ref{app:identities}. After a translation, we may also assume $x_\star=0$. By Algorithm~\ref{alg:tmmgd} and~\eqref{eq:output},
\[
\begin{aligned}
x_k&=x_{k-1}-\nabla f(x_{k-1})-m_k,&\qquad m_{k+1}&=\beta m_k+(\rho+\beta)\nabla f(x_k),\\
z_k&=x_{k-1}-\nabla f(x_{k-1})-\gamma_\star m_k,&\qquad z_{k+1}&=x_k-\nabla f(x_k)-\gamma_\star m_{k+1},
\end{aligned}
\]
and $z_1=x_0-\frac{1+q}{2q}\nabla f(x_0)$, since $m_1=\frac{\rho}{1-\rho}\nabla f(x_0)$ and $\gamma_\star\rho/(1-\rho)=(1-q)/(2q)$. With $\gamma_\star(\rho+\beta)=\rho/(1-\rho)$, $\gamma_\star(\rho-\beta)=\rho$ and $\sqrt q=1-\rho$, these expressions give
\[
\sqrt q\,(\rho z_k-z_{k+1})=g^\mu_k,
\qquad
\sqrt q\,(\rho z_k+z_{k+1})=g^L_k-\rho^2g^L_{k-1},
\qquad
2q\,z_1=qg^L_0-g^\mu_0.
\]
Hence $\rho^2\norm{z_k}^2-\norm{z_{k+1}}^2=\ip{g^\mu_k}{g^L_k-\rho^2g^L_{k-1}}/q$ and $4q^2\norm{z_1}^2=\norm{qg^L_0+g^\mu_0}^2-4q\ip{g^\mu_0}{g^L_0}$, where $qg^L_0+g^\mu_0=(1-q)\nabla f(x_0)$. Identity~\eqref{eq:twopoint}, with $(i,j)=(k-1,k)$ and $(k,k)$, and with $(i,j)=(0,0)$, writes these inner products with interpolation quantities, which gives~\eqref{eq:pot-decay} and~\eqref{eq:pot-start}. These steps and both identities can be checked symbolically, for any $L$ and $\mu$. The proof only uses the updates and~\eqref{eq:twopoint}, whose derivation only uses the zero gradient at $x_\star$, so both identities remain valid when $x_\star$ is replaced by any point that is assigned the gradient $0$, such as a virtual minimizer.
\end{proof}

\subsection{The hard function}
\label{app:hard}

The numbers $a_i$ and $\eta_i$ of Section~\ref{sec:tmmgd-minimax} are positive and satisfy, for $0\le i<N$,
\begin{equation}\label{eq:hard-recurrence}
q\,\eta_i<a_i\le\eta_i,
\qquad
a_{i+1}\eta_{i+1}=(\eta_i-a_i)(a_i-q\,\eta_i).
\end{equation}
The second identity follows by substitution: for $1\le i<N$, $\eta_i-a_i=\rho\,\eta_i$ and $a_i-q\eta_i=\sqrt q\,\rho\,\eta_i$, while for $i=0$ both sides equal $(1-q)^2/(4q)$.

Let $v_i:=a_ie_i-q\bar x_i$ for $0\le i\le N$, so that $\mathcal K$ is the convex hull of the points $(v_i-\eta)/(1-q)$. Since $a_N=q\eta_N$, $v_N=q\eta$.

\begin{proof}[Proof of Lemma~\ref{lem:hard}]
We use two results of \citet{drori2026weighted}, for $L=1$ and $\mu=q$. First, if $\mathcal K$ is the convex hull of finitely many points, then $F\in\mathcal F_{q,1}$, with gradient $\nabla F(x)=qx+\eta+(1-q)\Pi_{\mathcal K}(x)$, where $\Pi_{\mathcal K}$ is the projection onto $\mathcal K$, and $-\eta$ is a minimizer if $-\eta\in\mathcal K$, as here. Second, he builds such points from chain points, gradients at these points, and a minimizer. With $\bar x_i$, $a_ie_i$ and $-\eta$, his points are the points $(v_i-\eta)/(1-q)$, $0\le i\le N$, which is the case of his construction that recovers the construction of \citet{droritaylor2022oracle}. His conditions on the chain then reduce to $\eta_i>0$, $q\eta_i<a_i$ and $(\eta_i-a_i)(a_i-q\eta_i)\ge a_{i+1}\eta_{i+1}$ for $0\le i<N$, $0<a_N<\eta_N$, and $(\eta_{N-1}-a_{N-1})(a_{N-1}-q\eta_{N-1})\ge q\eta_N^2$, which hold by~\eqref{eq:hard-recurrence} and $a_N=q\eta_N$. Under these conditions, he proves that, for $0\le k<N$ and $x\in\operatorname{span}\{e_0,\dots,e_{k-1}\}$, $\nabla F(x)\in\operatorname{span}\{e_0,\dots,e_k\}$, which is the last claim of the lemma, and $\Pi_{\mathcal K}(x)$ lies in the convex hull $\mathcal K_k$ of the points $(v_i-\eta)/(1-q)$ with $i\le k$. He works in $\R^{N+1}$ and these properties carry over to $\R^d$, since $\mathcal K$ lies in $\operatorname{span}\{e_0,\dots,e_N\}$.

It remains to compute $\nabla F(\bar x_k)$ for $0\le k\le N$. By the above, $\Pi_{\mathcal K}(\bar x_k)$ is the projection of $\bar x_k$ onto $\mathcal K_k$, also for $k=N$, since $\mathcal K_N=\mathcal K$. The similarity $x\mapsto(1-q)x+\eta$ maps $\mathcal K_k$ onto the convex hull of $v_0,\dots,v_k$, and $\bar x_k$ to $(1-q)\bar x_k+\eta$, whose difference with $v_k$ is $(\eta_k-a_k)e_k+\sum_{t>k}\eta_te_t$. This difference has the inner product $-a_k(\eta_k-a_k)\le0$ with $v_j-v_k$ for $j<k$, so $v_k$ is the projection of $(1-q)\bar x_k+\eta$ onto this hull, and $\Pi_{\mathcal K}(\bar x_k)=(v_k-\eta)/(1-q)$. Hence $\nabla F(\bar x_k)=q\bar x_k+v_k=a_ke_k$.
\end{proof}

\section{Proofs for Section~\ref{sec:item}}
\label{app:item}

\subsection{The sequence \texorpdfstring{$A_k$}{A\_k}}
\label{app:item-identities}
 
Let $j\ge0$. The numerator in~\eqref{eq:item-A} is $\bigl(\sqrt{1+A_j}+\sqrt{1+qA_j}\bigr)^2$, which gives the first identity below. The second one follows, since its sides are positive and the difference of their squares is $q$ times that of the first identity:
\begin{equation}\label{eq:item-linear}
\begin{gathered}
(1-q)\sqrt{A_{j+1}}=\sqrt{1+A_j}+\sqrt{1+qA_j},\\
(1-q)\sqrt{1+qA_{j+1}}=\sqrt{1+qA_j}+q\sqrt{1+A_j}.
\end{gathered}
\end{equation}
Multiplying the first identity by $\sqrt{1+A_j}-\sqrt{1+qA_j}$ gives $\sqrt{A_{j+1}}\bigl(\sqrt{1+A_j}-\sqrt{1+qA_j}\bigr)=A_j$, and combining this relation with~\eqref{eq:item-linear} multiplied by $\sqrt{A_{j+1}}$ gives
\begin{equation}\label{eq:item-identities}
\begin{gathered}
(1+q)A_{j+1}-A_j=2\sqrt{A_{j+1}(1+qA_{j+1})},
\qquad
A_j+(1-q)A_{j+1}=2\sqrt{(1+A_j)A_{j+1}},\\
(1-q)A_{j+1}-A_j=2\sqrt{A_{j+1}(1+qA_j)}.
\end{gathered}
\end{equation}
These computations can be checked symbolically.
 
For the factor $R_N=\rho^{2N}(1+qA_N)/q$ of Remark~\ref{rem:price}: since $1-q=\rho(1+\sqrt q)$ and $q\sqrt{1+A_k}<\sqrt q\sqrt{1+qA_k}$, the second identity of~\eqref{eq:item-linear} gives $\rho\sqrt{1+qA_{k+1}}<\sqrt{1+qA_k}$, that is, $R_{k+1}<R_k$. Hence $R_N\le R_1=(1+q)^2/(q(1+\sqrt q)^2)<1/q$, since $A_1=4/(1-q)^2$. Adding $\sqrt q$ times the first identity of~\eqref{eq:item-linear} to the second gives $\rho\bigl(\sqrt{1+qA_{k+1}}+\sqrt{qA_{k+1}}\bigr)=\sqrt{1+qA_k}+\sqrt{q(1+A_k)}$, whose left-hand side is smaller than $\rho\bigl(\sqrt{1+qA_{k+1}}+\sqrt{q(1+A_{k+1})}\bigr)$. By induction from $A_0=0$, $\sqrt{1+qA_N}+\sqrt{q(1+A_N)}\ge(1+\sqrt q)/\rho^N$, and since $\sqrt{q(1+A_N)}<\sqrt{1+qA_N}$, it follows that $1+qA_N\ge(1+\sqrt q)^2/(4\rho^{2N})$, that is, $R_N\ge(1+\sqrt q)^2/(4q)>1/(4q)$.

\subsection{ITEM in momentum form}
\label{app:item-momentum}
 
In this subsection only, $y_k$, $z_k$, $\beta_k$ and $\delta_k$ denote the iterates and coefficients of \citet{taylordrori2023item}. Replacing $f$ by $f/L$ leaves the iterates unchanged and divides the momenta by $L$, so we may assume $L=1$. ITEM starts from $y_0=z_0=x_0$ and computes, for $k=0,\dots,N-1$,
\[
z_{k+1}=(1-q\delta_k)z_k+q\delta_ky_k-\delta_k\nabla f(y_k),
\qquad
y_{k+1}=(1-\beta_{k+1})z_{k+1}+\beta_{k+1}\bigl(y_k-\nabla f(y_k)\bigr),
\]
the second update only for $k\le N-2$, with $\beta_k=A_k/((1-q)A_{k+1})$ and $\delta_k=\sqrt{A_{k+1}/(1+qA_{k+1})}$, where we write $\delta_k$ as \citet{shuwang2026unified} do. It evaluates gradients at $y_0,\dots,y_{N-1}$ and returns $z_N$.
 
Let $x_k:=y_k$ for $k\le N-1$ and $x_N:=z_N$, let $\beta_N:=0$, and define $m_{k+1}$ by $x_{k+1}=x_k-\nabla f(x_k)-m_{k+1}$, with $m_0:=0$. Eliminating $z_k$ from the updates, with $z_k-x_k=-\beta_km_k/(1-\beta_k)$, gives $m_{k+1}=w_k\bigl(\nabla f(x_k)+\theta_km_k\bigr)$ with
\[
w_k=(1-\beta_{k+1})(\delta_k-1),
\qquad
\theta_k=\frac{(1-q\delta_k)\,\beta_k}{(\delta_k-1)(1-\beta_k)},
\]
and~\eqref{eq:item-linear} and~\eqref{eq:item-identities} give $\delta_k-1=\sqrt{1+A_k}/\sqrt{1+qA_{k+1}}$, $1-q\delta_k=\sqrt{1+qA_k}/\sqrt{1+qA_{k+1}}$ and $1-\beta_j=2\sqrt{1+qA_j}\big/\bigl((1-q)\sqrt{A_{j+1}}\bigr)$ for $0\le j\le N-1$. Inserting these expressions in $w_k$ and $\theta_k$ gives~\eqref{eq:item-coeffs}, where $\beta_N=0$ gives $w_{N-1}=\delta_{N-1}-1$. These computations can be checked symbolically. For $N\ge2$, the claim $w_0<w_{N-1}$ of Section~\ref{sec:itemgg} follows: $w_0=(1-\beta_1)(\delta_0-1)<\delta_0-1\le\delta_{N-1}-1=w_{N-1}$, since $0<\beta_1<1$, $\delta_0=2/(1+q)>1$ and $\delta_k$ increases with $k$, as $A_{k+1}$ does.

\subsection{Potential identities of ITEM-GG}
\label{app:itemgg-potential}
 
\begin{proof}[Proof of Lemma~\ref{lem:itemgg-potential}]
Replacing $f$ by $f/L$ leaves the iterates of Algorithm~\ref{alg:itemgg} unchanged, divides the momenta by $L$, and divides both sides of~\eqref{eq:itemgg-start} and~\eqref{eq:itemgg-decay} by $L^2$. We may therefore assume $L=1$ and $\mu=q$, and use the notation of Appendix~\ref{app:identities}. Algorithm~\ref{alg:itemgg} computes $m_1=w_{N-1}\nabla f(x_0)$, $m_{k+1}=w_{N-1-k}\bigl(\nabla f(x_k)+\theta_{N-k}m_k\bigr)$ for $k\ge1$, and $x_{k+1}=x_k-\nabla f(x_k)-m_{k+1}$. Let $\zeta^-_k:=\nabla f(x_k)+\theta_{N-k}m_k$ and $\zeta^+_k:=\nabla f(x_k)-(1-2\theta_{N-k})m_k$. Since $x_k^+-x_N^+=g^L_k-g^L_N$, the potential~\eqref{eq:itemgg-potential} reads
\[
\Phi_k=\norm{\zeta^-_k}^2+\ip{\zeta^+_k}{g^L_k-g^L_N}-(1-q)\,I_{N,k}.
\]
By~\eqref{eq:item-coeffs} and~\eqref{eq:item-identities}, the coefficients satisfy $w_n(1-\theta_n)=\sqrt{A_{n+1}/A_{n+2}}$ and $w_n(1+q-2\theta_n)=1-q-A_{n+1}/A_{n+2}$ for $0\le n\le N-2$, and $w_{N-1}(1-\theta_{N-1})=(1-q)\sqrt{A_N}\big/\bigl(2\sqrt{1+qA_N}\bigr)$ and $w_{N-1}(1+q-2\theta_{N-1})=1-q$. With these identities, the updates give, for $1\le k\le N-1$,
\[
\zeta^+_{k+1}-\frac{A_{N-k}}{A_{N-k+1}}\,\zeta^+_k=g^\mu_{k+1}-g^\mu_k,
\qquad
\ip{\zeta^+_{k+1}}{g^L_{k+1}-g^L_k}=\frac{A_{N-k}}{A_{N-k+1}}\norm{\zeta^-_k}^2-\norm{\zeta^-_{k+1}}^2,
\]
and, for the first step, $\zeta^+_1=g^\mu_1-g^\mu_0$ and $\norm{\zeta^-_1}^2+\ip{\zeta^+_1}{g^L_1-g^L_0}=\frac{(1-q)^2A_N}{4(1+qA_N)}\norm{\nabla f(x_0)}^2$. Inserting these identities into $\frac{A_{N-k}}{A_{N-k+1}}\Phi_k-\Phi_{k+1}$ and into $\Phi_1$ reduces~\eqref{eq:itemgg-decay} and~\eqref{eq:itemgg-start} to~\eqref{eq:threepoint}, with $(a,b,c)=(k,k+1,N)$ and $(0,1,N)$. These steps and both identities can be checked symbolically.
\end{proof}

\subsection{Rates of TMM-GG and of the output \texorpdfstring{$z^G_N$}{zGN}}
\label{app:momentum-bound}
\label{app:tmmgg-out}

Replacing $f$ by $f/L$ leaves the iterates of Algorithm~\ref{alg:tmmgd} and the ratios of Theorems~\ref{thm:tmmgg} and~\ref{thm:tmmgg-out} unchanged, so we may assume $L=1$ and $\mu=q$. Both rates of TMM-GG follow from the potential~\eqref{eq:potential} written against a virtual minimizer. Let $y$ be any point, with virtual minimizer $y^+=y-\nabla f(y)$. By Appendix~\ref{app:pot}, \eqref{eq:pot-start} and~\eqref{eq:pot-decay} hold with $x_\star$ replaced by $y^+$. On their right-hand sides, every interpolation quantity is between two iterates or has $y^+$ as the first point, and the latter are $I_{y^+,x_j}=I_{y,x_j}\ge0$ by~\eqref{eq:virtual}. As in the proof of Theorem~\ref{thm:rate}, it follows that, for every $j\ge1$,
\begin{equation}\label{eq:virtual-potential}
\norm{z_j-y^+}^2+\frac{1-q}{q}\,I_{x_{j-1},y^+}\le\rho^{2j-2}\Bigl(\frac{1-q}{2q}\Bigr)^2\norm{\nabla f(x_0)}^2.
\end{equation}
In the potential, $y^+$ is the second point of $I_{x_{j-1},y^+}$, so this term need not be nonnegative. Each proof below evaluates it.

\begin{proof}[Proof of Theorem~\ref{thm:tmmgg}]
We first bound the momenta. Fix $k\ge1$, and take $y=x_{k-1}$ and $j=k$ in~\eqref{eq:virtual-potential}. Then $I_{x_{k-1},x_{k-1}^+}=0$ by~\eqref{eq:virtual}, and $z_k-x_{k-1}^+=-\gamma_\star m_k$ by~\eqref{eq:output}. Since $(1-q)/(2q\gamma_\star)=\rho/\sqrt q$, it follows that $\norm{m_k}\le\rho^k\norm{\nabla f(x_0)}/\sqrt q$. For $f_\mu$ and $f_L$, every $I$ vanishes, so this bound is an equality.

Continuing the loop of Algorithm~\ref{alg:tmmgd} for one more step gives $m_{N+1}=\beta m_N+(\rho+\beta)\nabla f(x_N)$. Since $\beta/(\rho+\beta)=\rho/2$ and $1/(\rho+\beta)=(2-\rho)/(2\rho)$, it follows that $\nabla f(x_N)=\frac{2-\rho}{2\rho}\,m_{N+1}-\frac\rho2\,m_N$. By the bound on the momenta, applied to $m_{N+1}$ and $m_N$, the norm of the right-hand side is at most $\bigl(\frac{2-\rho}{2\rho}\rho^{N+1}+\frac\rho2\rho^N\bigr)\norm{\nabla f(x_0)}/\sqrt q=\rho^N\norm{\nabla f(x_0)}/\sqrt q$. For $f_L$, $\nabla f(x)=x-x_\star$ for every $x$, so line~4 with $k=N$ gives $\nabla f(x_N)=\nabla f(x_{N-1})-\nabla f(x_{N-1})-m_N=-m_N$, and the equality case of the bound on the momenta concludes.
\end{proof}

\begin{proof}[Proof of Theorem~\ref{thm:tmmgg-out}]
Take $y=z^G_N$ and $j=N$ in~\eqref{eq:virtual-potential}. By~\eqref{eq:output} and~\eqref{eq:output-gg}, $z_N-(z^G_N)^+=(\gamma_G-\gamma_\star)m_N+\nabla f(z^G_N)$ and $x_{N-1}^+-(z^G_N)^+=\gamma_Gm_N+\nabla f(z^G_N)$, so the second identity of~\eqref{eq:virtual} gives $I_{x_{N-1},(z^G_N)^+}=I_{x_{N-1},z^G_N}+\ip{\nabla f(z^G_N)}{\gamma_Gm_N+\nabla f(z^G_N)}$. The left-hand side of~\eqref{eq:virtual-potential} is therefore
\[
(\gamma_G-\gamma_\star)^2\norm{m_N}^2+\frac1q\norm{\nabla f(z^G_N)}^2+\frac{1-q}{q}\,I_{x_{N-1},z^G_N},
\]
where the terms in $\ip{m_N}{\nabla f(z^G_N)}$ cancel because $(1+q)\gamma_G=2q\gamma_\star$, and $\gamma_G$ is the only coefficient of the last momentum step for which they do. Since $I_{x_{N-1},z^G_N}\ge0$, it follows that $\norm{\nabla f(z^G_N)}^2\le\frac{(1-q)^2}{4q}\,\rho^{2N-2}\norm{\nabla f(x_0)}^2$, which is the claim, since $1-q=\rho(1+\sqrt q)$.
\end{proof}

\subsection{Minimax optimality and uniqueness of ITEM-GG}
\label{app:itemgg-minimax}

This appendix proves Theorem~\ref{thm:itemgg-minimax}. Replacing $f$ by $f/L$ leaves the iterates of a method~\eqref{eq:fsfom} and the ratio $\norm{\nabla f(x_N)}^2/\norm{\nabla f(x_0)}^2$ unchanged, so we assume $L=1$ and $\mu=q$, take $\bar L=1$, and use the notation and identities of Appendix~\ref{app:identities}. Here $N\ge1$, $0<q<1$, and $e_0,\dots,e_N$ is the canonical basis of $\R^{N+1}$. All matrices, vectors and scalar coefficients introduced below are local to this subsection. By~\eqref{eq:item-A}, $A_1=4/(1-q)^2$. The computations below can be checked symbolically, for small values of $N$ when they depend on $N$.

The proof eliminates the method from its Performance Estimation certificates, starting with the output row and then removing the remaining rows in reverse order. Each elimination bounds the change in a scalar terminal weight $\delta_m$. These bounds yield the sequence $A_k$, and their equality cases determine the multipliers and the method.

\emph{The method.} For $k=0,\dots,N$, let $\chi_k:=\nabla f(x_k)-(x_k-x_0)$. Since $x_k-x_0$ is a combination of $\nabla f(x_0),\dots,\nabla f(x_{k-1})$, there is a unit lower-triangular matrix $P$ with $\nabla f(x_k)=\sum_{i\le k}P_{ki}\chi_i$ for $k=0,\dots,N$. Conversely, every unit lower-triangular $P$ comes from exactly one $H$. Let $\pi:=P^\top e_N$ be the last row of $P$, which encodes the output step.

\emph{Certificates.} A certificate combines the inequalities $I_{ij}\ge0$, $0\le i\ne j\le N$, with weights $4\lambda_{ij}/(1-q)$, $\lambda_{ij}\ge0$, such that the function values cancel: $\sum_j\lambda_{ij}=\sum_j\lambda_{ji}$ for every $i$. Encode the multipliers in the matrix $\Lambda$ with $\Lambda_{kj}:=-\lambda_{jk}$ for $j\ne k$ and $\Lambda_{kk}:=\sum_j\lambda_{jk}$. Its off-diagonal entries are nonpositive, its rows sum to zero, and its columns sum to zero by the cancellation. Apply~\eqref{eq:twopoint} to every entry of $\Lambda$, diagonal included, where $I_{kk}=0$, and sum. The terms with $x_\star$ cancel, by the cancellation for $I_{k\star}$ and identically for $I_{\star k}$, and we get $(1-q)\sum_{i\ne j}\lambda_{ij}I_{ij}=\sum_{k,j}\Lambda_{kj}\ip{g^\mu_k}{g^L_j}$. Since the rows and columns of $\Lambda$ sum to zero, the right-hand side does not change if $x_\star$ is replaced by $x_0$ in $g^\mu$ and $g^L$, which adds $q(x_0-x_\star)$ to every $g^\mu_k$ and $x_\star-x_0$ to every $g^L_k$. Then $g^L_k=-\chi_k$ and $g^\mu_k=(1-q)\nabla f(x_k)+q\chi_k$. With $\sym(B):=\frac12(B+B^\top)$, $\nabla f(x_0)=\chi_0$ and $\nabla f(x_N)=\sum_i\pi_i\chi_i$, this gives, for every $\tau\in\R$,
\begin{equation}\label{eq:ggcert}
\frac1{A_1}\Bigl(\tau\norm{\nabla f(x_0)}^2-\norm{\nabla f(x_N)}^2\Bigr)-(1-q)\sum_{i\ne j}\lambda_{ij}I_{ij}=\sum_{k,j}\Psi_{kj}\ip{\chi_k}{\chi_j},
\end{equation}
\begin{equation}\label{eq:ggPsi}
\Psi:=\frac\tau{A_1}e_0e_0^\top-\frac1{A_1}\pi\pi^\top+(1-q)\sym(\Lambda^\top P)+q\sym(\Lambda).
\end{equation}
If $\Psi\succeq0$, then $\norm{\nabla f(x_N)}^2\le\tau\norm{\nabla f(x_0)}^2$ for every $f\in\mathcal F_{q,1}$ and every $x_0$, and we call $(\tau,\Lambda)$ a certificate for the method.

\emph{The output step.} For $v\in\R^{N+1}$, the output row enters $v^\top\Psi v$ only through the number $\pi^\top v=(Pv)_N$. Completing the square in this number, with $A_1(1-q)^2=4$, gives
\begin{equation}\label{eq:ggQ}
v^\top\Psi v=v^\top Qv-\frac1{A_1}\Bigl(\pi^\top v-\frac2{1-q}(\Lambda v)_N\Bigr)^2,
\end{equation}
where $Q$ is the symmetric matrix, independent of $\pi$, with
\[
v^\top Qv=\frac\tau{A_1}v_0^2+q\,v^\top\Lambda v+(1-q)\sum_{k<N}(\Lambda v)_k(Pv)_k+(\Lambda v)_N^2.
\]
For a certificate, $Q\succeq0$, and $\Lambda_{NN}\ge1/A_1$, since $\Psi_{NN}=\Lambda_{NN}-1/A_1$ by $\pi_N=1$ and the triangularity of $P$.

\begin{sublemma}[no duality gap, after \cite{taylorhendrickxglineur2017}]\label{lem:gg-duality}
Let $d\ge N+1$. The worst-case ratio of a method~\eqref{eq:fsfom} is the infimum of $\tau$ over its certificates, and also over its certificates with $\lambda_{ij}>0$ for all $i\ne j$.
\end{sublemma}

\begin{proof}
As in \citet{taylorhendrickxglineur2017}, who prove the absence of a duality gap for performance criteria under an initial distance condition when every $H_{k,k}\ne0$, Theorem~\ref{thm:interp} makes the worst-case ratio the value of a semidefinite program in the values $f(x_k)$ and the Gram matrix of the gradients, which has order $N+1$ because the ratio does not involve $x_\star$. Hence $d\ge N+1$ suffices. By~\eqref{eq:ggcert}, its dual is the minimization of $\tau$ over certificates. By conic duality, both claims follow from one certificate with $\Psi\succ0$ and positive multipliers: since $\Psi$ is affine in $(\tau,\Lambda)$, mixing a small amount of it into any certificate gives certificates with positive multipliers and nearly the same $\tau$. With $\epsilon_0:=\nabla f(x_0)$ and $\epsilon_{k+1}:=\nabla f(x_{k+1})-\nabla f(x_k)-\frac{1+q}2(x_{k+1}-x_k)$, a unit lower-triangular transformation of the gradients, $I_{k,k+1}+I_{k+1,k}=\frac{1-q}4\norm{x_{k+1}-x_k}^2-\frac1{1-q}\norm{\epsilon_{k+1}}^2$. So, for $\lambda_{k,k+1}=\lambda_{k+1,k}=\ell_k$, $0\le k<N$, $A_1$ times the left-hand side of~\eqref{eq:ggcert} is
\[
\tau\norm{\epsilon_0}^2-\norm{\nabla f(x_N)}^2+\sum_{k<N}\ell_k\bigl(A_1\norm{\epsilon_{k+1}}^2-\norm{x_{k+1}-x_k}^2\bigr).
\]
As $x_{k+1}-x_k$ only involves $\epsilon_0,\dots,\epsilon_k$, this form, and hence $\Psi$, is positive definite once $\ell_{N-1}$, then $\ell_{N-2},\dots,\ell_0$, and finally $\tau$ are large enough. Adding a small $\nu>0$ to every multiplier keeps $\Psi\succ0$ and the cancellation.
\end{proof}

\emph{Removing the remaining steps.} After completing the square in the output row, we eliminate rows $N-1,\dots,1$ of the method. At each step, a scalar Schur complement updates the multipliers, and restricting the quadratic form removes the corresponding row of $P$.

\begin{sublemma}[backward elimination]\label{lem:gg-reduction}
Assume that $\lambda_{ij}>0$ for all $i\ne j$ and $Q\succeq0$. Set $\mathcal I_m:=\{0,\dots,m,N\}$ for $0\le m\le N-1$ and $\Lambda^{(N-1)}:=\Lambda$. For $m=N-1,\dots,1$, define the matrix $\Lambda^{(m-1)}$ with indices in $\mathcal I_{m-1}$ by
\begin{equation}\label{eq:gg-elimination}
\Lambda^{(m-1)}_{ij}:=\Lambda^{(m)}_{ij}-\frac{\Lambda^{(m)}_{im}\Lambda^{(m)}_{mj}}{\Lambda^{(m)}_{mm}}\qquad(i,j\in\mathcal I_{m-1}).
\end{equation}
Every $\Lambda^{(m)}$ has strictly negative off-diagonal entries, positive diagonal entries, and zero row and column sums. The symmetric matrices $Q_m$ on coordinates $\mathcal I_m$, defined by
\begin{equation}\label{eq:gg-reducedQ}
v^\top Q_m v:=\frac\tau{A_1}v_0^2+q\,v^\top\Lambda^{(m)}v+(1-q)\sum_{k=0}^m(\Lambda^{(m)}v)_k\sum_{i=0}^kP_{ki}v_i+\bigl((\Lambda^{(m)}v)_N\bigr)^2,
\end{equation}
are positive semidefinite, with $Q_{N-1}=Q$. Write
\[
\alpha_m:=\Lambda^{(m)}_{mm},\qquad\beta_m:=-\Lambda^{(m)}_{mN},\qquad\gamma_m:=-\Lambda^{(m)}_{Nm},\qquad\delta_m:=\Lambda^{(m)}_{NN}.
\]
Then $0<\gamma_m\le\min(\alpha_m,\delta_m)$, $\beta_m>0$, $\alpha_0=\beta_0=\gamma_0=\delta_0$, and
\begin{equation}\label{eq:gg-scalar-update}
\delta_{N-1}=\Lambda_{NN},\qquad\delta_{m-1}=\delta_m-\frac{\beta_m\gamma_m}{\alpha_m}\quad(1\le m\le N-1).
\end{equation}
\end{sublemma}

\begin{proof}
For distinct $i,j\in\mathcal I_{m-1}$, \eqref{eq:gg-elimination} subtracts a positive number from $\Lambda^{(m)}_{ij}<0$, so $\Lambda^{(m-1)}_{ij}<0$. Summing~\eqref{eq:gg-elimination} over $j$, or over $i$, shows that the zero row and column sums are preserved. Thus every remaining diagonal entry, being the negative sum of the other entries in its row, is positive, which proves inductively that all divisions in~\eqref{eq:gg-elimination} are valid. To remove row $m$, restrict $Q_m$ by imposing $(\Lambda^{(m)}v)_m=0$, or equivalently
\begin{equation}\label{eq:gg-substitution}
v_m=-\frac{1}{\alpha_m}\sum_{i\in\mathcal I_{m-1}}\Lambda^{(m)}_{mi}v_i.
\end{equation}
On this subspace, the remaining coordinates of $\Lambda^{(m)}v$ equal $\Lambda^{(m-1)}v_{\mathcal I_{m-1}}$, so $v^\top\Lambda^{(m)}v$ is the quadratic form of $\Lambda^{(m-1)}$ at $v_{\mathcal I_{m-1}}$. The term containing row $m$ of $P$ vanishes, and the preceding rows do not involve $v_m$ because $P$ is lower triangular. Substituting in~\eqref{eq:gg-reducedQ} therefore gives exactly $Q_{m-1}$, which is positive semidefinite. Finally, the column sum at $m$ gives $\alpha_m-\gamma_m=-\sum_{i<m}\Lambda^{(m)}_{im}\ge0$, and the row sum at $N$ gives $\delta_m-\gamma_m=-\sum_{i<m}\Lambda^{(m)}_{Ni}\ge0$. For $m=0$, the balanced $2\times2$ matrix has all four coefficients equal. The scalar update is the $(N,N)$ entry of~\eqref{eq:gg-elimination}.
\end{proof}

The principal block of $Q_m$ on coordinates $m,N$ is independent of the method: when the preceding coordinates vanish, $\sum_{i=0}^kP_{ki}v_i$ is zero for $k<m$ and equals $v_m$ for $k=m$. Thus, for $1\le m\le N-1$,
\begin{equation}\label{eq:gg-corner}
B_m:=(Q_m)_{\{m,N\},\{m,N\}}=
\begin{pmatrix}
\alpha_m+\gamma_m^2&-\dfrac{\beta_m+q\gamma_m}{2}-\gamma_m\delta_m\\[3pt]
-\dfrac{\beta_m+q\gamma_m}{2}-\gamma_m\delta_m&\delta_m(q+\delta_m)
\end{pmatrix}\succeq0,
\end{equation}
and $B_0:=Q_0$ is given by the same formula with $m=0$ and $\tau/A_1$ added to its top-left entry.

\begin{sublemma}[one step]\label{lem:gg-onestep}
For $\delta>0$, let
\begin{equation}\label{eq:ggpsi}
\psi(\delta):=\frac{(1-q)^2}{\bigl(\sqrt{1+1/\delta}+\sqrt{1+q/\delta}\bigr)^2}=\delta\bigl(\sqrt{1+\delta}-\sqrt{q+\delta}\bigr)^2.
\end{equation}
Under the assumptions of Lemma~\ref{lem:gg-reduction}, the following hold:
\begin{enumerate}
\item[(a)] for $1\le m\le N-1$, with $u_m:=\gamma_m^2/\alpha_m\in[0,\delta_m]$,
\begin{equation}\label{eq:gg-square-bound}
\delta_{m-1}\ge\Bigl(\sqrt{\delta_m(1+u_m)}-\sqrt{(q+\delta_m)u_m}\Bigr)^2\ge\psi(\delta_m).
\end{equation}
For fixed $\delta_m$, the middle expression is strictly decreasing in $u_m\in[0,\delta_m]$. Equality throughout holds if and only if $\alpha_m=\gamma_m=\delta_m$ and $\beta_m=\delta_m-\psi(\delta_m)$;
\item[(b)] $\tau\ge\delta_0/(q+\delta_0)$;
\item[(c)] $\psi$ is strictly increasing, and $\psi(1/A_k)=1/A_{k+1}$ for every $k\ge1$.
\end{enumerate}
\end{sublemma}

\begin{proof}
(a) Omit the subscript $m$. Since $\det B_m\ge0$, equation~\eqref{eq:gg-corner} gives $(\beta+q\gamma)/2+\gamma\delta\le\sqrt{(\alpha+\gamma^2)\delta(q+\delta)}$. Multiplying by $2\gamma/\alpha$ and using $u=\gamma^2/\alpha$ and~\eqref{eq:gg-scalar-update} yields
\[
\delta_{m-1}=\delta-\frac{\beta\gamma}{\alpha}\ge\delta+(q+2\delta)u-2\sqrt{\delta(q+\delta)u(1+u)}=\Bigl(\sqrt{\delta(1+u)}-\sqrt{(q+\delta)u}\Bigr)^2.
\]
By Lemma~\ref{lem:gg-reduction}, $0\le u\le\gamma\le\delta$. On this interval,
\[
\sqrt{\delta(1+u)}-\sqrt{(q+\delta)u}=\frac{\delta-qu}{\sqrt{\delta(1+u)}+\sqrt{(q+\delta)u}}
\]
is positive and strictly decreasing: its numerator is positive and decreasing, and its denominator is positive and increasing. Its square is therefore minimized at $u=\delta$, where it equals $\psi(\delta)$. Equality forces $u=\gamma=\delta$, hence $\alpha=\delta$, and~\eqref{eq:gg-scalar-update} then gives $\beta=\delta-\psi(\delta)$. Conversely, these values give equality throughout, and $\det B_m=0$.

(b) Since $\alpha_0=\beta_0=\gamma_0=\delta_0>0$ and $A_1(1-q)^2=4$, the determinant of $B_0$ is $\delta_0\bigl(\tau(q+\delta_0)-\delta_0\bigr)/A_1\ge0$, which gives the claim, with $\det B_0=0$ at equality.

(c) Rationalizing the difference of square roots gives the two expressions in~\eqref{eq:ggpsi}. The denominator of the first is strictly decreasing in $\delta$, so $\psi$ is strictly increasing. For $\delta=1/A_k$, this denominator is $\bigl(\sqrt{1+A_k}+\sqrt{1+qA_k}\bigr)^2=(1-q)^2A_{k+1}$ by~\eqref{eq:item-linear}, so $\psi(1/A_k)=1/A_{k+1}$.
\end{proof}

\begin{proof}[Proof of the lower bound in Theorem~\ref{thm:itemgg-minimax}]
Let $(\tau,\Lambda)$ be a certificate with positive multipliers. The output step gives $\delta_{N-1}=\Lambda_{NN}\ge1/A_1$. Lemma~\ref{lem:gg-onestep}(a) and~(c) then give, successively for $m=N-1,\dots,1$, $\delta_{m-1}\ge\psi(\delta_m)\ge\psi(1/A_{N-m})=1/A_{N-m+1}$. Thus $\delta_0\ge1/A_N$, and Lemma~\ref{lem:gg-onestep}(b) yields
\[
\tau\ge\frac{\delta_0}{q+\delta_0}\ge\frac{1/A_N}{q+1/A_N}=\frac1{1+qA_N}.
\]
Lemma~\ref{lem:gg-duality} concludes. For $N=1$, no intermediate elimination is needed, and the same initial and final inequalities apply.
\end{proof}

\emph{The certificate of ITEM-GG.} Telescoping the identities of Lemma~\ref{lem:itemgg-potential} gives a certificate of ITEM-GG with $\tau=\bar\tau:=1/(1+qA_N)$ and $\Psi=0$. In its matrix $\Lambda^\star$, each column $k\le N-1$ has only two nonzero entries, $\Lambda^\star_{kk}=-\Lambda^\star_{k+1,k}=1/A_{N-k}$, and the last column follows from the zero row sums. By induction on $m$ from $N-1$ down to $0$, $\Lambda^{\star,(m)}$ consists of the rows $0,\dots,m$ of $\Lambda^\star$, restricted to the columns $\mathcal I_m$, and of a last row with entries $-1/A_{N-m}$ at $m$ and $1/A_{N-m}$ at $N$. Indeed, column $m$ of $\Lambda^{\star,(m)}$ has only two nonzero entries, $1/A_{N-m}$ and $-1/A_{N-m}$, the second one in row $N$, so~\eqref{eq:gg-elimination} adds row $m$ to row $N$ and deletes row and column $m$. Hence
\begin{align*}
\alpha_m^\star=\gamma_m^\star=\delta_m^\star&=\frac1{A_{N-m}}\quad(0\le m\le N-1),\\
\beta_0^\star=\frac1{A_N},\qquad\beta_m^\star&=\frac1{A_{N-m}}-\frac1{A_{N-m+1}}\quad(1\le m\le N-1).
\end{align*}
The pivots $\alpha_m^\star$ are positive, so the substitutions~\eqref{eq:gg-substitution}, the restrictions in the proof of Lemma~\ref{lem:gg-reduction} and the formula~\eqref{eq:gg-corner} remain valid despite the zero multipliers. All scalar bounds are equalities, and every $B_m$ is singular at $(\bar\tau,\Lambda^\star)$. These blocks do not depend on $P$.

\begin{proof}[Proof of the uniqueness in Theorem~\ref{thm:itemgg-minimax}]
Let a method with matrix $P$ attain $\bar\tau$. By Lemma~\ref{lem:gg-duality}, it has certificates $(\tau_j,\Lambda_j)$ with positive multipliers and $\tau_j\to\bar\tau$. We first show that $\Lambda_j\to\Lambda^\star$, and then that $(\bar\tau,\Lambda^\star)$ certifies only ITEM-GG. We do not assume that an optimal certificate exists.

\emph{The scalar bounds become equalities.} Use a superscript $j$ for the scalar coefficients obtained from $\Lambda_j$. The lower-bound proof gives $\delta_m^j\ge1/A_{N-m}$ and $\tau_j\ge\delta_0^j/(q+\delta_0^j)$, so $\delta_0^j\to1/A_N$. Applying $\psi(\delta_m^j)\le\delta_{m-1}^j$ successively and using the strict monotonicity and continuity of $\psi$ gives $\delta_m^j\to1/A_{N-m}$ for every $m$. For $m\ge1$, both outer terms in~\eqref{eq:gg-square-bound} converge to $1/A_{N-m+1}$. Since the middle expression is continuous and has its unique minimum at $u_m=\delta_m$, it follows that $u_m^j\to1/A_{N-m}$: the sequence $u_m^j\in[0,\delta_m^j]$ is bounded, and any other subsequential limit would give a strict inequality. Now $u_m^j\le\gamma_m^j\le\delta_m^j$, $\alpha_m^j=(\gamma_m^j)^2/u_m^j$ and $\beta_m^j=\alpha_m^j(\delta_m^j-\delta_{m-1}^j)/\gamma_m^j$ imply $(\alpha_m^j,\beta_m^j,\gamma_m^j,\delta_m^j)\to(\alpha_m^\star,\beta_m^\star,\gamma_m^\star,\delta_m^\star)$. For $m=0$, all four coefficients equal $\delta_0^j\to1/A_N$.

\emph{Reversing the elimination recovers the multipliers.} Let $\Lambda_j^{(m)}$ denote the reduced matrices. We prove by induction on $m$ that $\Lambda_j^{(m)}\to\Lambda^{\star,(m)}$. For $m=0$, this follows from $\Lambda_j^{(0)}=\delta_0^j\left(\begin{smallmatrix}1&-1\\-1&1\end{smallmatrix}\right)$. For $m\ge1$, the entries of $\Lambda_j^{(m)}$ above the pivot in column $m$ and before column $m$ in row $N$ are nonpositive, and the zero sums of column $m$ and row $N$ show that they add up to $\gamma_m^j-\alpha_m^j\to0$ and $\gamma_m^j-\delta_m^j\to0$, respectively. Hence they tend to zero. For $i<m$, \eqref{eq:gg-elimination} gives
\[
(\Lambda_j^{(m-1)})_{Ni}=(\Lambda_j^{(m)})_{Ni}+\frac{\gamma_m^j}{\alpha_m^j}(\Lambda_j^{(m)})_{mi}.
\]
By the induction hypothesis and $\gamma_m^j/\alpha_m^j\to1$, this determines the limit of row $m$. The corrections in~\eqref{eq:gg-elimination} to the rows $i<m$ then tend to zero, which determines their limits as well. All limits are therefore uniquely determined by $\Lambda^{\star,(m-1)}$ and the limiting scalar coefficients. The same relations hold for $\Lambda^{\star,(m)}$, whose entries above the pivot in column $m$ and before column $m$ in row $N$ are zero. This proves the induction, and in particular $\Lambda_j\to\Lambda^\star$. By continuity, $(\bar\tau,\Lambda^\star)$ is a certificate for $P$.

\emph{The pair $(\bar\tau,\Lambda^\star)$ certifies only ITEM-GG.} Let $P$ be a unit lower-triangular matrix certified by this pair, with output row $\pi$, and let $P^\star$ and $\pi^\star$ be those of ITEM-GG. For $0\le m\le N-1$, the singular block $B_m$ has the positive bottom-right entry $\delta_m^\star(q+\delta_m^\star)$, so $B_m(1,t_m)^\top=0$ for some $t_m\in\R$. Let $w_m\in\R^{N+1}$ have the entries $0$ before $m$, $1$ at $m$ and $t_m$ at $N$, and the entries $m+1,\dots,N-1$ given successively by the substitutions~\eqref{eq:gg-substitution}. The proof of Lemma~\ref{lem:gg-reduction} then gives $w_m^\top Qw_m=(1,t_m)B_m(1,t_m)^\top=0$, so~\eqref{eq:ggQ} and $\Psi\succeq0$ yield
\[
\Psi w_m=0,\qquad\pi^\top w_m=\frac2{1-q}(\Lambda^\star w_m)_N.
\]
Moreover, $\Psi e_N=0$, since $\Psi_{NN}=\Lambda^\star_{NN}-1/A_1=0$. The vectors $w_0,\dots,w_{N-1},e_N$ are the columns of a unit lower-triangular matrix, hence a basis. Thus $\Psi=0$. Moreover, $\pi$ is determined by $\pi_N=1$ and the values $\pi^\top w_m$, which the display gives independently of $P$. As ITEM-GG is also certified by this pair, $\pi=\pi^\star$. Now $\Delta:=P-P^\star$ is strictly lower triangular with zero last row, and subtracting the two zero certificate matrices gives $\sym(\Lambda^{\star\top}\Delta)=0$. Set $X:=\Lambda^{\star\top}\Delta$. For $0\le k<N$,
\[
X_{kl}=\frac{\Delta_{kl}-\Delta_{k+1,l}}{A_{N-k}}\qquad(0\le l\le N).
\]
Both entries on the right vanish when $l>k$. Thus $X$ is lower triangular and skew-symmetric, hence zero. The display then shows that consecutive rows of $\Delta$ are equal. Since its first row is zero, $\Delta=0$ and $P=P^\star$.
\end{proof}

\bibliography{Limited_memory_paper}

@article{bansalgupta2019,
  author  = {Bansal, N. and Gupta, A.},
  title   = {Potential-function proofs for gradient methods},
  journal = {Theory Comput.},
  volume  = {15},
  number  = {4},
  pages   = {1--32},
  year    = {2019}
}

@article{barretaylorbach2023,
  author  = {Barr{\'e}, M. and Taylor, A. B. and Bach, F.},
  title   = {Principled analyses and design of first-order methods with inexact proximal operators},
  journal = {Math. Program.},
  volume  = {201},
  number  = {1--2},
  pages   = {185--230},
  year    = {2023}
}

@article{bokaltschuler2026,
  author  = {Bok, J. and Altschuler, J. M.},
  title   = {Optimized methods for composite optimization: a reduction perspective},
  journal = {Math. Program.},
  year    = {2026},
  doi     = {10.1007/s10107-026-02377-7}
}

@article{dasguptavanparysryu2024,
  author  = {Das Gupta, S. and Van Parys, B. P. G. and Ryu, E. K.},
  title   = {Branch-and-bound performance estimation programming: a unified methodology for constructing optimal optimization methods},
  journal = {Math. Program.},
  volume  = {204},
  number  = {1--2},
  pages   = {567--639},
  year    = {2024}
}

@book{dewildeveen1998,
  author    = {Dewilde, P. and {van der Veen}, A.-J.},
  title     = {Time-Varying Systems and Computations},
  publisher = {Kluwer},
  address   = {Boston},
  year      = {1998}
}

@article{diakonikolaswang2022,
  author  = {Diakonikolas, J. and Wang, P.},
  title   = {Potential function-based framework for minimizing gradients in convex and min-max optimization},
  journal = {SIAM J. Optim.},
  volume  = {32},
  pages   = {1668--1697},
  year    = {2022}
}

@article{drori2017exact,
  author  = {Drori, Y.},
  title   = {The exact information-based complexity of smooth convex minimization},
  journal = {J. Complexity},
  volume  = {39},
  pages   = {1--16},
  year    = {2017}
}

@misc{drori2026weighted,
  author       = {Drori, Y.},
  title        = {Exact oracle complexity for function-value suboptimality under weighted initial conditions},
  howpublished = {arXiv:2609.33360v2},
  year         = {2026}
}

@article{droritaylor2020efficient,
  author  = {Drori, Y. and Taylor, A. B.},
  title   = {Efficient first-order methods for convex minimization: a constructive approach},
  journal = {Math. Program.},
  volume  = {184},
  number  = {1--2},
  pages   = {183--220},
  year    = {2020}
}

@article{droritaylor2022oracle,
  author  = {Drori, Y. and Taylor, A. B.},
  title   = {On the oracle complexity of smooth strongly convex minimization},
  journal = {J. Complexity},
  volume  = {68},
  pages   = {101590},
  year    = {2022}
}

@article{droriteboulle2014,
  author  = {Drori, Y. and Teboulle, M.},
  title   = {Performance of first-order methods for smooth convex minimization: a novel approach},
  journal = {Math. Program.},
  volume  = {145},
  number  = {1--2},
  pages   = {451--482},
  year    = {2014}
}

@article{eidelmangohberg1999,
  author  = {Eidelman, Y. and Gohberg, I.},
  title   = {On a new class of structured matrices},
  journal = {Integral Equations Operator Theory},
  volume  = {34},
  number  = {3},
  pages   = {293--324},
  year    = {1999}
}

@book{eidelmangohberghaimovici2014,
  author    = {Eidelman, Y. and Gohberg, I. and Haimovici, I.},
  title     = {Separable Type Representations of Matrices and Fast Algorithms. Volume~1: Basics. Completion Problems. Multiplication and Inversion Algorithms},
  series    = {Oper. Theory Adv. Appl.},
  volume    = {234},
  publisher = {Birkh{\"a}user},
  address   = {Basel},
  year      = {2014}
}

@article{fazlyabribeiromoraripreciado2018,
  author  = {Fazlyab, M. and Ribeiro, A. and Morari, M. and Preciado, V. M.},
  title   = {Analysis of optimization algorithms via integral quadratic constraints: nonstrongly convex problems},
  journal = {SIAM J. Optim.},
  volume  = {28},
  number  = {3},
  pages   = {2654--2689},
  year    = {2018}
}

@article{goujaud2024pepit,
  author  = {Goujaud, B. and Moucer, C. and Glineur, F. and Hendrickx, J. M. and Taylor, A. B. and Dieuleveut, A.},
  title   = {{PEPit}: computer-assisted worst-case analyses of first-order optimization methods in {P}ython},
  journal = {Math. Program. Comput.},
  volume  = {16},
  number  = {3},
  pages   = {337--367},
  year    = {2024}
}

@misc{grimmerjopark2026,
  author       = {Grimmer, B. and Jo, S. and Park, C.},
  title        = {On minimax optimality and uniqueness of fixed-step first-order methods for smooth convex optimization},
  howpublished = {arXiv:2609.07671},
  year         = {2026}
}

@incollection{kaashoekwoerdeman1989,
  author    = {Kaashoek, M. A. and Woerdeman, H. J.},
  title     = {Minimal lower separable representations: characterization and construction},
  editor    = {Dym, H. and Goldberg, S. and Kaashoek, M. A. and Lancaster, P.},
  booktitle = {The Gohberg Anniversary Collection, Vol.~II: Topics in Analysis and Operator Theory},
  series    = {Oper. Theory Adv. Appl.},
  volume    = {41},
  pages     = {329--344},
  publisher = {Birkh{\"a}user},
  address   = {Basel},
  year      = {1989}
}

@misc{kamrihendrickxglineur2025,
  author       = {Kamri, Y. and Hendrickx, J. M. and Glineur, F.},
  title        = {Numerical design of optimized first-order algorithms},
  howpublished = {arXiv:2507.20773},
  year         = {2025}
}

@article{kimfessler2016ogm,
  author  = {Kim, D. and Fessler, J. A.},
  title   = {Optimized first-order methods for smooth convex minimization},
  journal = {Math. Program.},
  volume  = {159},
  number  = {1--2},
  pages   = {81--107},
  year    = {2016}
}

@article{kim2021accelerated,
  author  = {Kim, D.},
  title   = {Accelerated proximal point method for maximally monotone operators},
  journal = {Math. Program.},
  volume  = {190},
  number  = {1--2},
  pages   = {57--87},
  year    = {2021}
}

@article{kimfessler2021ogmg,
  author  = {Kim, D. and Fessler, J. A.},
  title   = {Optimizing the efficiency of first-order methods for decreasing the gradient of smooth convex functions},
  journal = {J. Optim. Theory Appl.},
  volume  = {188},
  number  = {1},
  pages   = {192--219},
  year    = {2021}
}

@misc{kimryudasgupta2026,
  author       = {Kim, H. and Ryu, E. K. and Das Gupta, S.},
  title        = {A domain-specific harness for end-to-end automation of optimization research},
  howpublished = {arXiv:2608.07407},
  year         = {2026}
}

@inproceedings{kim2023timereversed,
  author    = {Kim, J. and Ozdaglar, A. E. and Park, C. and Ryu, E. K.},
  title     = {Time-reversed dissipation induces duality between minimizing gradient norm and function value},
  booktitle = {Advances in Neural Information Processing Systems 36},
  pages     = {23389--23440},
  year      = {2023}
}

@inproceedings{kornowski2024open,
  author    = {Kornowski, G. and Shamir, O.},
  title     = {Open problem: anytime convergence rate of gradient descent},
  booktitle = {Proc. 37th Conf. Learning Theory},
  series    = {PMLR},
  volume    = {247},
  pages     = {5335--5339},
  year      = {2024}
}

@inproceedings{leeparkryu2021,
  author    = {Lee, J. and Park, C. and Ryu, E. K.},
  title     = {A geometric structure of acceleration and its role in making gradients small fast},
  booktitle = {Advances in Neural Information Processing Systems 34},
  pages     = {11999--12012},
  year      = {2021}
}

@article{lessardrechtpackard2016,
  author  = {Lessard, L. and Recht, B. and Packard, A.},
  title   = {Analysis and design of optimization algorithms via integral quadratic constraints},
  journal = {SIAM J. Optim.},
  volume  = {26},
  number  = {1},
  pages   = {57--95},
  year    = {2016}
}

@inproceedings{lessardseiler2020,
  author    = {Lessard, L. and Seiler, P.},
  title     = {Direct synthesis of iterative algorithms with bounds on achievable worst-case convergence rate},
  booktitle = {Proc. American Control Conf.},
  pages     = {119--125},
  year      = {2020}
}

@article{lessardudell2026,
  author  = {Lessard, L. and Udell, M.},
  title   = {Algebraic characterization of equivalence between oracle-based iterative algorithms},
  journal = {Math. Program.},
  year    = {2026},
  doi     = {10.1007/s10107-026-02402-9}
}

@book{nemirovskyyudin1983,
  author    = {Nemirovsky, A. S. and Yudin, D. B.},
  title     = {Problem Complexity and Method Efficiency in Optimization},
  publisher = {Wiley},
  address   = {Chichester},
  year      = {1983}
}

@article{nesterov1983,
  author  = {Nesterov, Y. E.},
  title   = {A method of solving a convex programming problem with convergence rate {$O(1/k^2)$}},
  journal = {Soviet Math. Dokl.},
  volume  = {27},
  number  = {2},
  pages   = {372--376},
  year    = {1983}
}

@book{nesterov2004,
  author    = {Nesterov, Y.},
  title     = {Introductory Lectures on Convex Optimization: A Basic Course},
  series    = {Appl. Optim.},
  volume    = {87},
  publisher = {Kluwer},
  address   = {Boston},
  year      = {2004}
}

@article{odonoghuecandes2015,
  author  = {O'Donoghue, B. and Cand{\`e}s, E.},
  title   = {Adaptive restart for accelerated gradient schemes},
  journal = {Found. Comput. Math.},
  volume  = {15},
  number  = {3},
  pages   = {715--732},
  year    = {2015}
}

@article{rodmanwoerdeman1995,
  author  = {Rodman, L. and Woerdeman, H. J.},
  title   = {Perturbations, singular values, and ranks of partial triangular matrices},
  journal = {SIAM J. Matrix Anal. Appl.},
  volume  = {16},
  number  = {1},
  pages   = {278--288},
  year    = {1995}
}

@article{schererebenbauer2021,
  author  = {Scherer, C. W. and Ebenbauer, C.},
  title   = {Convex synthesis of accelerated gradient algorithms},
  journal = {SIAM J. Control Optim.},
  volume  = {59},
  number  = {6},
  pages   = {4615--4645},
  year    = {2021}
}

@misc{shuwang2026unified,
  author       = {Shu, K. and Wang, A. L.},
  title        = {A unified theory of {$H$}-duality in first-order methods},
  howpublished = {arXiv:2609.03281},
  year         = {2026}
}

@article{taylordrori2023item,
  author  = {Taylor, A. B. and Drori, Y.},
  title   = {An optimal gradient method for smooth strongly convex minimization},
  journal = {Math. Program.},
  volume  = {199},
  number  = {1--2},
  pages   = {557--594},
  year    = {2023}
}

@article{taylorhendrickxglineur2017,
  author  = {Taylor, A. B. and Hendrickx, J. M. and Glineur, F.},
  title   = {Smooth strongly convex interpolation and exact worst-case performance of first-order methods},
  journal = {Math. Program.},
  volume  = {161},
  number  = {1--2},
  pages   = {307--345},
  year    = {2017}
}

@inproceedings{taylorvanscoylessard2018,
  author    = {Taylor, A. B. and Van Scoy, B. and Lessard, L.},
  title     = {{L}yapunov functions for first-order methods: tight automated convergence guarantees},
  booktitle = {Proc. 35th ICML},
  series    = {PMLR},
  volume    = {80},
  pages     = {4897--4906},
  year      = {2018}
}

@misc{upadhyaya2026proxitem,
  author       = {Upadhyaya, M. and Berg Thomsen, D. and Dieuleveut, A. and Taylor, A. B.},
  title        = {An optimal first-order method for smooth and strongly convex composite optimization and its stationary limit},
  howpublished = {arXiv:2605.22929},
  year         = {2026}
}

@article{vandebrilvanbarelmastronardi2005,
  author  = {Vandebril, R. and Van Barel, M. and Mastronardi, N.},
  title   = {A note on the representation and definition of semiseparable matrices},
  journal = {Numer. Linear Algebra Appl.},
  volume  = {12},
  number  = {8},
  pages   = {839--858},
  year    = {2005}
}

@book{vandebrilvanbarelmastronardi2008,
  author    = {Vandebril, R. and Van Barel, M. and Mastronardi, N.},
  title     = {Matrix Computations and Semiseparable Matrices, Vol.~1: Linear Systems},
  publisher = {Johns Hopkins University Press},
  address   = {Baltimore},
  year      = {2008}
}

@article{vanderveendewilde1994lowcomplexity,
  author  = {{van der Veen}, A.-J. and Dewilde, P.},
  title   = {On low-complexity approximation of matrices},
  journal = {Linear Algebra Appl.},
  volume  = {205--206},
  pages   = {1145--1201},
  year    = {1994}
}

@article{vanscoy2018tmm,
  author  = {Van Scoy, B. and Freeman, R. A. and Lynch, K. M.},
  title   = {The fastest known globally convergent first-order method for minimizing strongly convex functions},
  journal = {IEEE Control Syst. Lett.},
  volume  = {2},
  number  = {1},
  pages   = {49--54},
  year    = {2018}
}

@misc{vernimmen2026anytime,
  author       = {Vernimmen, P. and Glineur, F.},
  title        = {Can acceleration in gradient-norm minimization be anytime? {S}harp last-iterate limits in smooth convex optimization},
  howpublished = {arXiv:2609.05710},
  year         = {2026}
}

@inproceedings{weibel2026optimized,
  author    = {Weibel, J. and Gaillard, P. and Koolen, W. M. and Taylor, A. B.},
  title     = {Optimized projection-free algorithms for online learning: construction and worst-case analysis},
  booktitle = {Proc. 29th Int. Conf. Artificial Intelligence and Statistics},
  series    = {PMLR},
  volume    = {300},
  pages     = {5122--5130},
  year      = {2026}
}

@article{woerdeman1987,
  author  = {Woerdeman, H. J.},
  title   = {The lower order of lower triangular operators and minimal rank extensions},
  journal = {Integral Equations Operator Theory},
  volume  = {10},
  number  = {6},
  pages   = {859--879},
  year    = {1987}
}

@inproceedings{yoonkimsuhryu2024,
  author    = {Yoon, T. and Kim, J. and Suh, J. J. and Ryu, E. K.},
  title     = {Optimal acceleration for minimax and fixed-point problems is not unique},
  booktitle = {Proc. 41st ICML},
  series    = {PMLR},
  volume    = {235},
  pages     = {57244--57314},
  year      = {2024}
}

@misc{yoonryugrimmer2025,
  author       = {Yoon, T. and Ryu, E. K. and Grimmer, B.},
  title        = {{H}-invariance theory: a complete characterization of minimax optimal fixed-point algorithms},
  howpublished = {arXiv:2511.14915},
  year         = {2025}
}

@misc{yoongrimmer2026,
  author       = {Yoon, T. and Grimmer, B.},
  title        = {A theory of composition and duality of extremal optimal fixed-point algorithms},
  howpublished = {arXiv:2605.02231},
  year         = {2026}
}

@inproceedings{zhang2025anytime,
  author    = {Zhang, Z. and Lee, J. D. and Du, S. S. and Chen, Y.},
  title     = {Anytime acceleration of gradient descent},
  booktitle = {Proc. 38th Conf. Learning Theory},
  series    = {PMLR},
  volume    = {291},
  pages     = {5991--6013},
  year      = {2025}
}

\end{document}